\documentclass[10pt,a4paper]{article}
\usepackage{bbm}

\usepackage[leqno]{amsmath}
\usepackage{amsfonts}
\usepackage{graphicx}
\usepackage{amsmath}
\usepackage{amssymb}
\usepackage{latexsym}
\usepackage{amsmath, amsfonts,amssymb, amsthm, euscript,makeidx,color,mathrsfs}
\usepackage{enumerate}

\usepackage[colorlinks,linkcolor=blue,anchorcolor=green,citecolor=red]{hyperref}

\def\sqr#1#2{{\vcenter{\vbox{\hrule height.#2pt
              \hbox{\vrule width.#2pt height#1pt \kern#1pt \vrule width.#2pt}
              \hrule height.#2pt}}}}
\def\signed #1{{\unskip\nobreak\hfil\penalty50
              \hskip2em\hbox{}\nobreak\hfil#1
              \parfillskip=0pt \finalhyphendemerits=0 \par}}
\def\endpf{\signed {$\sqr69$}}

\def\5n{\negthinspace \negthinspace \negthinspace \negthinspace \negthinspace }
\def\4n{\negthinspace \negthinspace \negthinspace \negthinspace }
\def\3n{\negthinspace \negthinspace \negthinspace }
\def\2n{\negthinspace \negthinspace }
\def\1n{\negthinspace }

\def\dbE{\mathbb{E}}
\def\dbF{\mathbb{F}}

\def\dbK{\mathbb{K}}
\def\dbL{\mathbb{L}}

\def\dbP{\mathbb{P}}

\def\dbR{\mathbb{R}}

\def\dbU{\mathbb{U}}

\def\dbX{\mathbb{X}}
\def\dbY{\mathbb{Y}}

\def\={\buildrel \triangle \over =}

\def\ds{\displaystyle}

\def\ns{\noalign{\ss}}
\def\a{\alpha}
\def\b{\beta}
\def\g{\gamma}
\def\d{\delta}
\def\e{\varepsilon}
\def\z{\zeta}

\def\m{\mu}
\def\n{\nu}
\def\si{\sigma}
\def\t{\tau}
\def\f{\varphi}
\def\th{\theta}
\def\o{\omega}

\def\G{\Gamma}

\def\Th{\Theta}
\def\L{\Lambda}

\def\F{\Phi}
\def\O{\Omega}
\def\cF{{\cal F}}

\def\cN{{\cal N}}

\def\cR{{\cal R}}

\def\cU{{\cal U}}

\def\cY{{\cal Y}}

\def\no{\noindent}

\def\ss{\smallskip}
\def\ms{\medskip}

\def\q{\quad}
\def\qq{\qquad}
\def\hb{\hbox}

\def\ua{\mathop{\uparrow}}
\def\da{\mathop{\downarrow}}
\def\Ra{\mathop{\Rightarrow}}

\def\lan{\langle}
\def\ran{\rangle}

\def\h{\widehat}
\def\wt{\widetilde}

\def\cd{\cdot}
\def\cds{\cdots}

\def\ae{\hbox{\rm a.e.}}
\def\as{\hbox{\rm a.s.}}

\def\tr{\hbox{\rm tr$\,$}}

\def\les{\leqslant}
\def\ges{\geqslant}

\def\({\Big (}
\def\){\Big )}
\def\[{\Big[}
\def\]{\Big]}

\def\bde{\begin{definition}\label}
\def\ede{\end{definition}}
\def\be{\begin{equation}}
\def\bel{\begin{equation}\label}
\def\ee{\end{equation}}
\def\bt{\begin{theorem}\label}
\def\et{\end{theorem}}
\def\bc{\begin{corollary}\label}
\def\ec{\end{corollary}}
\def\bl{\begin{lemma}\label}
\def\el{\end{lemma}}
\def\bp{\begin{proposition}\label}
\def\ep{\end{proposition}}
\def\bas{\begin{assumption}\label}
\def\eas{\end{assumption}}
\def\br{\begin{remark}\label}
\def\er{\end{remark}}
\def\bex{\begin{example}\label}
\def\ex{\end{example}}
\def\ba{\begin{array}}
\def\ea{\end{array}}
\def\ed{\end{document}}

\def\square#1{\vbox{\hrule\hbox{\vrule height#1%
     \kern#1\vrule}\hrule}}
\def\rectangle#1#2{\vbox{\hrule\hbox{\vrule height#1%
     \kern#2\vrule}\hrule}}

\font\tenbb=msbm10 \font\sevenbb=msbm7 \font\fivebb=msbm5

\newfam\bbfam
\scriptscriptfont\bbfam=\fivebb \textfont\bbfam=\tenbb
\scriptfont\bbfam=\sevenbb

\newtheorem{theorem}{\hskip 1.3em Theorem}[section]
\newtheorem{definition}[theorem]{\hskip 1.3em Definition}
\newtheorem{proposition}[theorem]{\hskip 1.3em Proposition}
\newtheorem{corollary}[theorem]{\hskip 1.3em Corollary}
\newtheorem{lemma}[theorem]{\hskip 1.3em Lemma}
\newtheorem{remark}[theorem]{\hskip 1.3em Remark}
\newtheorem{example}[theorem]{\hskip 1.3em Example}

\newtheorem{assumption}[theorem]{\hskip 1.3em Assumption}

\makeatletter
   
   \@addtoreset{equation}{section}
\makeatother

\begin{document}

\title{\bf Exact Controllability of Linear Stochastic Differential Equations and Related Problems\thanks{This work is supported in part by
the National Natural Science Foundation of China (11471192,
11371375, 11526167), the Fundamental Research Funds for the Central
Universities (SWU113038, XDJK2014C076), the Nature Science
Foundation of Shandong Province (JQ201401), the Natural Science
Foundation of CQCSTC (2015jcyjA00017), and NSF Grant DMS-1406776.}}

\author{Yanqing Wang\footnote{School of Mathematics and Statistics, Southwest University, Chongqing 400715, China; email: {\tt yqwang@amss.ac.cn}},~~
Donghui Yang\footnote{School of Mathematics and Statistics, Central
South University, Changsha 410075, China; email: {\tt
donghyang@139.com}},~~Jiongmin Yong\footnote{Department of
Mathematics, University of Central Florida, Orlando, FL 32816, USA;
email: {\tt jiongmin.yong@ucf.edu}},~~and~~Zhiyong
Yu\footnote{Corresponding author, School of Mathematics, Shandong
University, Jinan 250100, China; email: {\tt yuzhiyong@sdu.edu.cn}}}

\maketitle

\begin{abstract}
A notion of $L^p$-exact controllability is introduced for linear
controlled (forward) stochastic differential equations, for which
several sufficient conditions are established. Further, it is proved
that the $L^p$-exact controllability, the validity of an
observability inequality for the adjoint equation, the solvability
of an optimization problem, and the solvability of an $L^p$-type norm
optimal control problem are all equivalent.

\end{abstract}

\ms

\no\bf Keywords: \rm Controlled stochastic differential equation,
$L^p$-exact controllability, observability inequality, norm optimal control problem.

\ms

\no\bf AMS subject classification: \rm 93B05, 93E20, 60H10

\maketitle

\section{Introduction}

Let $(\O,\cF,\dbF,\dbP)$ be a complete filtered probability space on
which a $d$-dimensional standard Brownian motion
$W(\cd)\equiv(W_1(\cd),\cds,W_d(\cd))$ is defined such that
$\dbF\equiv\{\cF_t\}_{t\ges0}$ is its natural filtration augmented
by all the $\dbP$-null sets. Consider the following linear
controlled (forward) stochastic differential equation (FSDE, for
short):
\bel{FSDE1}dX(t)=\[A(t)X(t)+B(t)u(t)\]dt+\sum_{k=1}^d\[C_k(t)X(t)+D_k(t)u(t)\]dW_k(t),\qq t\ges0,\ee
where $A,C_1,\cds,C_d:[0,T]\times\O\to\dbR^{n\times n}$ and
$B:D_1,\cds,D_d:[0,T]\times\O\to\dbR^{n\times m}$ are suitable
matrix-valued processes. Here, $\dbR^{n\times n}$ and $\dbR^{n\times
m}$ are the sets of all $(n\times n)$ and $(n\times m)$ real
matrices, respectively. In the above, $X(\cd)$ is the {\it state
process} valued in the $n$-dimensional (real) Euclidean space
$\dbR^n$ and $u(\cd)$ is the {\it control process} valued in the
$m$-dimensional (real) Euclidean space $\dbR^m$. We will denote
system (\ref{FSDE1}) by $[A(\cd),C(\cd);B(\cd),D(\cd)]$, with
$C(\cd)=(C_1(\cd),\cds,C_d(\cd))$ and $D(\cd)=(D_1(\cd),\cds,
D_d(\cd))$, and denote
$$[A(\cd),0;B(\cd),0]=[A(\cd);B(\cd)].$$
In addition, if $A(\cd)$ and $B(\cd)$ are deterministic, $[A(\cd);B(\cd)]$ is reduced to a linear controlled ordinary differential equation (ODE, for short), for which a very mature theory is available; See, for example, Wonham \cite{Wonham 1985}, and references cited therein.

\ms

For system \eqref{FSDE1}, a control process $u:[0,T]\times\O\to\dbR^m$ is said to be {\it feasible} on $[0,T]$ if $t\mapsto u(t)$ is $\dbF$-progressively measurable, and for any initial state $x\in\dbR^n$, the state equation \eqref{FSDE1} admits a unique strong solution $X(\cd)\equiv X(\cd\,;x,u(\cd))$ on $[0,T]$ with the property
\bel{|X|}X(0)=x,\qq\dbE\[\sup_{t\in[0,T]}|X(t)|\]<\infty.\ee
Let $\cU[0,T]$ be the set of all {\it feasible controls} on $[0,T]$, whose precise definition will be given a little later. Let $\dbU[0,T]\subseteq\cU[0,T]$. System $[A(\cd),C(\cd);B(\cd),D(\cd)]$ is said to be $L^p$-{\it exactly controllable} on $[0,T]$ by $\dbU[0,T]$ if for any initial state $x\in\dbR^n$ and any {\it terminal state} $\xi\in L^p_{\cF_T}(\O;\dbR^n)$, the set of all $\cF_T$-measurable random variables $\xi:\O\to\dbR^n$ with $\dbE|\xi|^p<\infty$, there exists a $u(\cd)\in\dbU[0,T]$ such that
\bel{X=xi}X(T;x,u(\cd))=\xi.\ee
Clearly, the above notion is a natural extension of controllability
for ODE system $[A(\cd);B(\cd)]$ (\cite{Wonham 1985}).

\ms

Controllability is one of the most important concepts in control
theory. For time-invariant linear ODE systems, it is well-known that
the controllability is equivalent to many conditions, among which,
the Kalman's rank condition is the most interesting one. It is also
known that the controllability of a controlled linear ODE system is
equivalent to the observability of its adjoint system. For
controlled linear partial differential equations (PDEs, for short),
the notion of controllability can further be split into the so-called {\it exact
controllability}, {\it null controllability}, {\it approximate
controllability}, which are not equivalent in general, and all of
these three are closely related to the so-called {\it unique
continuation} property, and/or the {\it observability inequality}
for the adjoint equation. For extensive surveys of controllability
results on deterministic systems, see Lee--Markus \cite{Lee-Markus
1967} for ODE systems, and Russell \cite{Russell
1978}, Lions \cite{Lions 1988a, Lions 1988b} and Zuazua \cite{Zuazua
2006} for PDE systems.

\ms

The study of the controllability for stochastic systems can be traced back to the work of Connor \cite{Connors-1967} in 1967, followed by Sunahara--Aihara--Kishino \cite{Sunahara-Aihara-Kishino 1975}, Zabczyk \cite{Zabczyk 1981}, Ehrhardt--Kliemann \cite{Ehrhardt-Kliemann-1982}, and Chen--Li--Peng--Yong \cite{Chen-Li-Peng-Yong-1993}. With the help of backward stochastic differential equations (BSDEs, for short), Peng \cite{Peng 1994} introduced the so-called {\it exact terminal controllability} and {\it exact controllability}\footnote{In the terminology of the current paper, this should roughly be called the $L^2$-exact controllability.} for linear FSDE system with constant coefficients; the former was characterized by the non-degeneracy of the matrix $D$ and the latter was characterized by a generalized version of Kalman's rank condition. Later, Liu--Peng \cite{Liu-Peng 2010} extended the results to the linear FSDE system with bounded time-varying coefficients, using a random version of Gramian matrix. On the other hand, Buckdahn--Quincampoix--Tessitore \cite{Buckdahn-Quincampoix-Tessitore 2006} and Goreac \cite{Goreac 2008} studied the so-called {\it approximate controllability} (in $L^2$ sense) for linear FSDE systems, with constant coefficients and with degenerate matrix $D$ in the state equation; Some generalized Kalman type conditions are obtained to characterize the approximate controllability. In \cite{Lu-Yong-Zhang 2012}, L\"u--Yong--Zhang established a representation of It\^o's integral as a Lebesgue/Bochner integral, which has some interesting consequences on controllability of a linear FSDE system with vanished diffusion (see below).

\ms

In this paper, for any $p\in [1,\infty)$, we propose a notion of
$L^p$-exact controllability (see Definition \ref{Sec3_Def_Lp_Exact})
for FSDE systems. When $p=2$ and all the coefficients of the
system are bounded, our notion is reduced to the one studied in
\cite{Peng 1994, Liu-Peng 2010}. We point out that since the
coefficients $B(\cd)$ and $D_k(\cd)$ are allowed to be unbounded,
the corresponding set of admissible controls is delicate and it
makes the controllability problem under consideration more
interesting (see below for detailed explanation). Inspired by the
results of deterministic systems, for any $p\in(1,\infty)$, we
introduce a stochastic version of observability inequality (see
Theorem \ref{equivalence}) for the adjoint equation, the validity of
which is proved to be equivalent to the $L^p$-exact controllability
of the linear FSDE (with random coefficients). This provides an
approach to study the controllability of linear FSDE systems
by establishing an inequality for BSDEs. Moreover, we introduce a family of
optimization problems of the adjoint linear BSDE (see Problem (O)
and Problem (O)$'$ in Section \ref{S_Observability}), and the
solvability of these optimization problems is proved to be
equivalent to the $L^p$-exact controllability of the linear FSDE. In
other words, we additionally provide an approach to study the exact
controllability through infinite-dimensional optimization theory.

\ms

As an application, we consider some $L^p$-type norm optimal control
problems (see Problem (N) and Problem (N)$'$ in Section
\ref{S_Norm_Optimal_Control}). The norm optimal control problem for
deterministic finite or infinite dimensional systems has been
investigated by many researchers (see e.g. \cite{Fattorini 1999,
Fattorini 2011, Gugat-Leugering 2008, Wang-Xu 2013, Wang-Zuazua
2012}). Recently, Gashi \cite{Gashi 2015} studied a norm optimal
control problem (in $L^2$ sense) for linear FSDE systems with
deterministic time-varying coefficients by virtue of the
corresponding Hamiltonian system and Riccati equation. Moreover,
Wang--Zhang \cite{Wang-Zhang 2015} studied a kind of approximately
norm optimal control problems for linear FSDEs. In the present
paper, with the help of optimization problems for BSDE, we solve the
norm optimal control problem for linear FSDE systems with random
coefficients (see Theorems \ref{Sec4_Theorem_Equiv} and
\ref{Sec4_weak_Theorem_Equiv}, and Corollary \ref{Corollary 5.5}).

\ms

The rest of this paper is organized as follows. In Section
\ref{S_Preliminary}, we present some preliminaries. Section \ref{Controllability} is devoted to the introduction of the $L^p$-exact controllability for linear FSDE systems. Some sufficient conditions of the $L^p$-exact controllability are established for two types of
systems: The diffusion is control-free and the diffusion is ``fully'' controlled.
In Section \ref{S_Observability}, we establish the equivalence among the
$L^p$-exact controllability, the
validity of observability inequality for the adjoint equation and
the solvability of optimal control problems. Finally, as an application, a norm optimization problem is considered in Section \ref{S_Norm_Optimal_Control}.

\section{Preliminaries}\label{S_Preliminary}

Recall that $\dbR^n$ is the $n$-dimensional (real) Euclidean
(vector) space with the standard Euclidean norm $|\cd|$ induced by the
standard Euclidean inner product $\lan\cd\,,\cd\ran$, and
$\dbR^{n\times m}$ is the space of all $(n\times m)$ (real)
matrices, with the inner product
$$\lan A,B\ran=\tr[A^\top B],\qq\forall A,B\in\dbR^{n\times m},$$
so that $\dbR^{n\times m}$ is also a Euclidean space. Hereafter, the
superscript $^\top$ denotes the transpose of a vector or a matrix. We
now introduce some spaces, besides $L^p_{\cF_T}(\O;\dbR^n)$
introduced in the previous section. Let $H=\dbR^n,\dbR^{n\times m}$,
etc., and $p,q\in[1,\infty)$.
\begin{itemize}
\item $L^p_\dbF(\O;L^q(0,T;H))$ is the set of all $\dbF$-progressively measurable processes $\f(\cd)$ valued in $H$ such that
$$\|\f(\cd)\|_{L^p_\dbF(\O;L^q(0,T;H))}\equiv\[\dbE\(\int_0^T|\f(t)|^qdt\)^{p\over q}\]^{1\over p}<\infty.$$
\item $L^q_\dbF(0,T;L^p(\O;H))$ is the set of all $\dbF$-progressively measurable processes $\f(\cd)$ valued in $H$ such that
$$\|\f(\cd)\|_{L^q_\dbF(0,T;L^p(\O;H))}\equiv\[\int_0^T\(\dbE|\f(t)|^p\)^{q\over p}dt\]^{1\over q}<\infty.$$
\item $L^p_\dbF(\O;C([0,T];H)$ is the set of all $\dbF$-progressively measurable processes $\f(\cd)$ valued in $H$ such that for almost all $\o\in\O$, $t\mapsto\f(t,\o)$ is continuous and
$$\|\f(\cd)\|_{L^p_\dbF(\O;C([0,T];H))}\equiv\[\dbE\(\sup_{t\in[0,T]}|X(t)|^p\)\]^{1\over p}<\infty.$$
\end{itemize}
In the similar manner, one may define
$$\ba{ll}
\ns\ds L^p_\dbF(\O;L^\infty(0,T;H)),\q
L^\infty_\dbF(0,T;L^p(\O;H)),\q
L^p_\dbF(0,T;L^\infty(\O;H)),\q L^\infty_\dbF(\O;L^p(0,T;H)),\\
\ns\ds L^\infty_\dbF(\O;L^\infty(0,T;H)),\q
L^\infty_\dbF(0,T;L^\infty(\O;H)),\q
L^\infty_\dbF(\O;C([0,T];H)).\ea$$
We have
\bel{L=L}L^p_\dbF(\O;L^p(0,T;H))=L^p_\dbF(0,T;L^p(\O;H))\equiv
L^p_\dbF(0,T;H),\qq p\in[1,\infty],\ee
and
\bel{Minkowski}\left\{\2n\ba{ll}
\ns\ds L^p_\dbF(\O;L^q(0,T;H))\subseteq L^q_\dbF(0,T;L^p(\O;H)),\qq1\les p\les q\les\infty,\\
\ns\ds L^q_\dbF(0,T;L^p(\O;H))\subseteq
L^p_\dbF(\O;L^q(0,T;H)),\qq1\les q\les p\les\infty.\ea\right.\ee
In particular,
\bel{L(1p)}L^1_\dbF(0,T;L^p(\O;H))\subseteq
L^p_\dbF(\O;L^1(0,T;H)),\qq1\les p\les\infty.\ee
Also, we have that
\bel{C-ne-C}\ba{ll}
\ns\ds L^p_\dbF(\O;C([0,T];H))\subseteq
L^p_\dbF(\O;L^\infty(0,T;H))\subseteq L^\infty_\dbF(0,T;L^p(\O;H)).\ea\ee
In fact, for $1\les p\les q<\infty$, by Minkowski's integral
inequality, we have
$$\ba{ll}
\ns\ds\|\f(\cd)\|_{L^q_\dbF(0,T;L^p(\O;H))}^p=\[\int_0^T\(\dbE|\f(t)|^p\)^{q\over
p}dt\]^{p\over q}\les\dbE\(\int_0^T|\f(t)|^qdt\)^{p\over
q}=\|\f(\cd)\|_{L^p_\dbF(\O;L^q(0,T;H))}^p.\ea$$
This gives the first inclusion in \eqref{Minkowski}. Other cases can
be proved similarly. 
Now, we introduce the following definition.

\begin{definition}\label{feasible control}
\rm An $\dbF$-progressive measurable process
$u:[0,T]\times\O\to\dbR^m$ is called a {\it feasible control} of
system $[A(\cd),C(\cd);B(\cd),D(\cd)]$ if under $u(\cd)$, for any
$x\in\dbR^n$, system $[A(\cd),C(\cd);B(\cd),D(\cd)]$ admits a unique
strong solution $X(\cd)\in L^1_\dbF(\O;C([0,T];\dbR^n))$ satisfying
$X(0)=x$. The set of feasible controls is denoted by $\cU[0,T]$.
\end{definition}

Now, for the state equation \eqref{FSDE1}, we introduce the
following basic hypothesis.

\ms

\bf (H1) \rm The $\dbR^{n\times n}$-valued processes
$A(\cd),C_1(\cd),\cds,C_d(\cd)$ satisfy
\bel{AC}A(\cd),C_1(\cd),\cds,C_d(\cd)\in
L^\infty_\dbF(0,T;\dbR^{n\times n}).\ee
The $\dbR^{n\times m}$-valued processes
$B(\cd),D_1(\cd),\cds,D_d(\cd)$ are $\dbF$-progressively measurable.

\ms

The following result gives a big class of feasible controls for
system $[A(\cd),C(\cd);B(\cd),D(\cd)]$, whose proof is standard.

\begin{proposition}\label{well-posedness of SDE}
\sl Let {\rm(H1)} hold. Let $u:[0,T]\times\O\to\dbR^n$ be
$\dbF$-progressively measurable. Suppose the following holds:
\bel{Sec2_Ups}p\ges 1,~B(\cd)u(\cd)\in
L^p_\dbF(\O;L^1(0,T;\dbR^n)),~D_k(\cd)u(\cd)\in
L^p_\dbF(\O;L^2(0,T;\dbR^n)),~1\les k\les d.\ee
Then $u(\cd)\in\cU[0,T]$, and the solution $X(\cd)\equiv
X(\cd\,;x,u(\cd))$ of \eqref{FSDE1} with initial state $x$ under control $u(\cd)$ satisfies the following:
\bel{Lp-estimate}
\|X(\cd)\|_{L^p_\dbF(\O;C([0,T];\dbR^n))}\1n \les\1n K\Big\{|x|\1n
+\1n\|B(\cd)u(\cd)\|_{L^p_\dbF(\O;L^1(0,T;\dbR^n))}\1n+\1n\sum_{k=1}^d\|D_k(\cd)u(\cd)\|_{L^p_\dbF
(\O;L^2(0,T;\dbR^n))}\Big\}.\ee
Hereafter, $K>0$ will denote a generic constant which could be
different from line to line. Further, if
\bel{Sec2_Upw}p\ges 2,~B(\cd)u(\cd)\in
L^1_\dbF(0,T;L^p(\O;\dbR^n)),~D_k(\cd)u(\cd)\in
L^2_\dbF(0,T;L^p(\O;\dbR^n)),~1\les k\les d,\ee
then $u(\cd)\in\cU[0,T]$, and the following holds:
\begin{equation}\label{Lp-estimate*}
\|X(\cd)\|_{L^p_\dbF(\O;C([0,T];\dbR^n))}\1n\les K\1n\Big\{|x|\1n
+\1n\|B(\cd)u(\cd)\|_{L^1_\dbF(0,T;L^p(\O;\dbR^n))}\1n+\1n\sum_{k=1}^d\|D_k(\cd)u(\cd)\|_{L^2_\dbF
(0,T;L^p(\O;\dbR^n))}\Big\}.
\end{equation}
\end{proposition}

\ms

The above result leads us to the following definitions.

\begin{definition}\label{Sec2_Lp feasible control}
\rm A control $u(\cdot)\in\cU[0,T]$ is said to be {\it $L^p$-feasible} (respectively, {\it $L^p$-restricted feasible}) for
system $[A(\cdot),C(\cdot);B(\cdot),D(\cdot)]$ if \eqref{Sec2_Ups}
(respectively, \eqref{Sec2_Upw}) holds true. The set of
$L^p$-feasible controls (respectively, $L^p$-restricted feasible
controls) is denoted by $\mathcal U^p[0,T]$ (respectively, $\mathcal
U^p_r[0,T]$).
\end{definition}

Now, let us introduce the following two sets of hypotheses.

\ms

{\bf(H2)} For some $\rho\in(1,\infty]$ and $\si\in(2,\infty]$, the
following hold:
\bel{BD}\ba{ll}
\ns\ds B(\cd)\in\left\{\2n\ba{ll}
\ns\ds L^\rho_\dbF(\O;L^{2\si\over \si+2}(0,T;\dbR^{n\times
m})),\q\rho\in
(1,\infty],\ \si\in (2,\infty),\\
\ns\ds L^\rho_\dbF(\O;L^2(0,T;\dbR^{n\times
m})),\q\rho\in(1,\infty],\ \si=\infty, \ea\right.\\ [3mm]
\ns\ds D_1(\cd),\cds,D_d(\cd)\in
L^\rho_\dbF(\O;L^\si(0,T;\dbR^{n\times m})).\ea\ee

\ms

{\bf(H2)$'$} For some $\rho,\si\in(2,\infty]$, the following hold:
\bel{BD'}\ba{ll}
\ns\ds B(\cd)\in\left\{\2n\ba{ll}
\ns\ds L^{2\si\over \si+2}_\dbF(0,T;L^\rho(\O;\dbR^{n\times
m})),\q\rho\in
(2,\infty],\ \si\in (2,\infty),\\
\ns\ds L^2_\dbF(0,T;L^\rho(\O;\dbR^{n\times
m})),\q\rho\in(2,\infty],\ \si=\infty, \ea\right.\\ [3mm]
\ns\ds D_1(\cd),\cds,D_d(\cd)\in
L^\si_\dbF(0,T;L^\rho(\O;\dbR^{n\times m})).\ea\ee

\ms

In what follows, (H2) will be used for problems involving
$L^p$-feasible controls and (H2)$'$ will be used for problems
involving $L^p$-restricted feasible controls. We now denote the
following set of controls:
$$\dbU^{p,\rho,\si}[0,T]=\left\{\2n\ba{ll}
\ns\ds L^{\rho p\over\rho-p}_\dbF(\O;L^{2\si\over \si-2}(0,T;\dbR^m)),\qq p\in[1,\rho),~\rho\in(1,\infty),~\si\in(2,\infty),\\
\ns\ds L^p_\dbF(\O;L^{2\si\over \si-2}(0,T;\dbR^m)),\qq p\in[1,\rho),~\rho=\infty,~\si\in(2,\infty),\\
\ns\ds L^{\rho p\over\rho-p}_\dbF(\O;L^2(0,T;\dbR^m)),\qq p\in[1,\rho),~\rho\in(1,\infty),~\si=\infty,\\
\ns\ds
L^p_\dbF(\O;L^2(0,T;\dbR^m)),\qq p\in[1,\rho),~\rho=\si=\infty.\ea\right.$$
Clearly, as $\rho\ua+\infty$ and/or $\si\ua+\infty$, the set
$\dbU^{p,\rho,\si}[0,T]$ is getting larger and larger. Whereas, as
$p\ua\rho$, the set $\dbU^{p,\rho,\si}[0,T]$ is getting smaller and
smaller. Similarly, we introduce
$$\dbU_r^{p,\rho,\si}[0,T]=\left\{\2n\ba{ll}
\ns\ds L^{2\si\over \si-2}_\dbF(0,T;L^{\rho p\over\rho-p}(\O;\dbR^m)),\qq p\in[2,\rho),~\rho,\si\in(2,\infty),\\
\ns\ds L^{2\si\over \si-2}_\dbF(0,T;L^p(\O;\dbR^m)),\qq p\in[2,\rho),~\rho=\infty,~\si\in(2,\infty),\\
\ns\ds L^2_\dbF(0,T;L^{\rho p\over\rho-p}(\O;\dbR^m)),\qq p\in[2,\rho),~\rho\in(2,\infty),~\si=\infty,\\
\ns\ds L^2_\dbF(0,T;L^p(\O;\dbR^m)),\qq p\in[2,\rho),~\rho=\si=\infty.\ea\right.$$
We have the following proposition.

\bp{}\sl {\rm(i)} Let {\rm(H1)} and {\rm(H2)} hold. Then
\bel{U in U}\dbU^{p,\rho,\si}[0,T]\subseteq\cU^p[0,T], \qq\forall
p\in[1,\rho).\ee

{\rm(ii)} Let {\rm(H1)} and {\rm(H2)$'$} hold. Then
\bel{U in U*}\dbU_r^{p,\rho,\si}[0,T]\subseteq\cU_r^p[0,T],
\qq\forall p\in[2,\rho).\ee\ep

\it Proof. \rm (i) Let (H1) and (H2) hold. Let
$u(\cd)\in\dbU^{p,\rho,\si}[0,T]$. We only consider the case that
$\rho<\infty$ and $\si<\infty$; Others can be proved similarly.
We make the following calculations (noting $\si>2$ and $1\les
p<\rho$),
\bel{|Bu|}\ba{ll}
\ns\ds\|B(\cd)u(\cd)\|_{L^p_\dbF(\O;L^1(0,T;\dbR^n))}^p=\dbE\(\int_0^T|B(t)u(t)|dt\)^p\\
\ns\ds\qq\qq\qq\qq\qq\q\les\dbE\[\(\int_0^T|B(t)|^{2\si\over \si+2}dt\)^{(\si+2)p\over2\si}\(\int_0^T|u(t)|^{2\si\over \si-2}dt\)^{(\si-2)p\over 2\si}\]\\
\ns\ds\qq\qq\qq\qq\qq\q\les\[\dbE\(\int_0^T|B(t)|^{2\si\over \si+2}dt\)^{(\si+2)\rho\over2\si}\]^{p\over\rho}\[\dbE\(\int_0^T|u(t)|^{2\si\over \si-2}dt\)^{(\si-2)p\rho\over 2\si(\rho-p)}\]^{\rho-p\over\rho}\\
\ns\ds\qq\qq\qq\qq\qq\q\equiv\|B(\cd)\|_{L^\rho_\dbF(\O;L^{2\si\over
\si+2}(0,T;\dbR^{n\times m}))}^p
\|u(\cd)\|_{L^{p\rho\over\rho-p}_\dbF(\O;L^{2\si\over
\si-2}(0,T;\dbR^m))}^p,\ea\ee
and for each $k=1,\cds,d$,
\bel{|Du|}\ba{ll}
\ns\ds\|D_k(\cd)u(\cd)\|_{L^p_\dbF(\O;L^2(0,T;\dbR^n))}^p=\dbE\(\int_0^T|D_k(t)u(t)|^2dt\)^{p\over2}\\
\ns\ds\qq\qq\qq\qq\qq\qq\les\dbE\[\(\int_0^T|D_k(t)|^\si dt\)^{p\over \si}\(\int_0^T|u(t)|^{2\si\over \si-2}dt\)^{(\si-2)p\over 2\si}\]\\
\ns\ds\qq\qq\qq\qq\qq\qq\les\[\dbE\(\int_0^T|D_k(t)|^\si dt\)^{\rho\over \si}\]^{p\over\rho}\[\dbE\(\int_0^T|u(t)|^{2\si\over \si-2}dt\)^{(\si-2)p\rho\over 2\si(\rho-p)}\]^{\rho-p\over\rho}\\
\ns\ds\qq\qq\qq\qq\qq\qq\equiv\|D_k(\cd)\|_{L^\rho_\dbF(\O;L^\si(0,T;\dbR^{n\times
m}))}^p \|u(\cd)\|_{L^{p\rho\over\rho-p}_\dbF(\O;L^{2\si\over
\si-2}(0,T;\dbR^m))}^p.\ea\ee
By Definition \ref{Sec2_Lp feasible control}, $u(\cd)\in\cU^p[0,T]$,
proving (i).

\ms

In the similar manner, we are able to prove (ii).
\endpf

\section{Exact Controllability}\label{Controllability}

We now give a precise definition of $L^p$-exact controllability.

\bde{Sec3_Def_Lp_Exact} \rm Let $\dbU[0,T]\subseteq\cU[0,T]$. System
$[A(\cd),C(\cd);B(\cd),D(\cd)]$ is said to be $L^p$-{\it exactly
controllable} by $\dbU[0,T]$ on the time interval $[0,T]$, if for
any $(x,\xi)\in\dbR^n\times L^p_{\cF_T}(\O;\dbR^n)$, there exists a
$u(\cd)\in\dbU[0,T]$ such that the solution $X(\cd)\in L^1_{\mathbb
F}(\Omega;C([0,T];\mathbb R^n))$ of \eqref{FSDE1} with $X(0)=x$
satisfies $X(T)=\xi$.

\ede

In the above, $\dbU[0,T]$ could be $\cU^q[0,T]$, or $\cU^q_r[0,T]$
for some suitable $q\ges1$, and also it could be
$\dbU^{p,\rho,\si}[0,T]$, or $\dbU^{p,\rho,\si}_r[0,T]$. We
emphasize that in defining the system to be $L^p$-exactly
controllable (for $p\ges1$) by $\dbU[0,T]$, we only require
$X(\cd)\in L^1_{\mathbb F}(\Omega;C([0,T];\mathbb R^n))$ (since
$\dbU[0,T]\subseteq\cU[0,T]$). Depending on the choice of
$\dbU[0,T]$, $X(\cd)$ might have better integrability/regularity but
we do not require any better property than $L^1_{\mathbb
F}(\Omega;C([0,T];\mathbb R^n))$. We will see shortly that this
gives us a great reflexibility.

\subsection{The case $D(\cd)=0$}

In this subsection, we consider system $[A(\cd),C(\cd);B(\cd),0]$,
i.e., the state equation reads
\bel{FSDE0}dX(t)=\[A(t)X(t)+B(t)u(t)\]dt+\sum_{k=1}^dC_k(t)X(t)dW_k(t),\qq
t\ges0.\ee
Thus, the control $u(\cd)$ does not appear in the diffusion. When
all the coefficients in the above are constants, it was shown in
\cite{Buckdahn-Quincampoix-Tessitore 2006} that the system is
approximately controllable (under some additional conditions) in the
following sense: For any $(x,\xi)\in\dbR^n\times
L^2_{\cF_T}(\O;\dbR^n)$, and any $\e>0$, there exists a
$u_\e(\cd)\in L^2_\dbF(0,T;\dbR^m)\equiv\dbU^{2,\infty,\infty}[0,T]$
such that the solution $X(\cd)\equiv X(\cd\,;x,u_\e(\cd))$ with
$X(0)=x$ satisfies
$$\|X(T)-\xi\|_{L^2_{\cF_T}(\O;\dbR^n)}<\e.$$
The following is our first result which improves the results of
\cite{Buckdahn-Quincampoix-Tessitore 2006} significantly.

\bt{D=0}\sl Let $D(\cdot)=0$ and {\rm(H1)} hold. Let
%
%
%
\bel{B(t)B(t)>d}B(t)B(t)^\top\ges\d I>0,\qq t\in[0,T],~\as,\ee
for some $\d>0$. Then for any $p>1$, system
$[A(\cd),C(\cd);B(\cd),0]$ is $L^p$-exactly controllable on $[0,T]$
by $\cU^q[0,T]$ with any $q\in(1,p)$. \et

\it Proof. \rm Consider the following system:
$$\left\{\2n\ba{ll}
\ns\ds dX(t)=v(t)dt+\sum_{k=1}^dC_k(t)X(t)dW_k(t),\qq t\ges0,\\
\ns\ds X(0)=x,\ea\right.$$
with $v(\cd)\in L^q_\dbF(\O;L^1(0,T;\dbR^n))$, $q>1$. Then the
unique solution $X(\cd)$ satisfies
$$\dbE\[\sup_{t\in[0,T]}|X(t)|^q\]\les K\[|x|^q+\dbE\(\int_0^T|v(t)|dt\)^q\].$$
Let $\F(\cd)$ be the solution to the following:
$$\left\{\2n\ba{ll}
\ns\ds d\F(t)=\sum_{k=1}^dC_k(t)\F(t)dW_k(t),\qq t\ges0,\\
\ns\ds\F(0)=I.\ea\right.$$
Then $\F(\cd)^{-1}$ exists and satisfies the following:
$$\left\{\2n\ba{ll}
\ns\ds d\big[\F(t)^{-1}\big]=\sum_{k=1}^d\F(t)^{-1}C_k(t)^2dt-\sum_{k=1}^d\F(t)^{-1}C_k(t)dW_k(t),\qq t\ges0,\\
\ns\ds\F(0)^{-1}=I.\ea\right.$$
Therefore, for any $q\ges1$,
$$\dbE\[\sup_{[t\in[0,T]}|\F(t)|^q+\sup_{t\in[0,T]}|\F(t)^{-1}|^q\]\les K(T,q),$$
with the constant $K(T,q)$ depending on $T$ and $q$ (as well as the
bound of $C_k(\cd)$), and we have the following variation of
constants formula for $X(\cd)$:
\bel{X(t)}X(t)=\F(t)x+\F(t)\int_0^t\F(s)^{-1}v(s)ds,\qq t\ges0.\ee
Now, for any $q\in(1,p)$, we want to choose some $v(\cd)\in
L^q_\dbF(\O;L^1(0,T;\dbR^n))$ so that $X(T)=\xi$ which is equivalent
to the following:
$$\F(T)^{-1}\xi-x=\int_0^T\h v(s)ds,$$
with
$$\h v(t)=\F(t)^{-1}v(t),\qq t\in[0,T].$$
Since $\xi\in L^p_{\cF_T}(\O;\dbR^n)$, for any $\bar q\in(q,p)$, we
have
$$\dbE|\F(T)^{-1}\xi|^{\bar q}\les\(\dbE|\F(T)^{-1}|^{\bar qp\over p-\bar q}\)^{p-\bar q\over p}\(\dbE|\xi|^p\)^{\bar q\over p}\les K\(\dbE|\xi|^p\)^{\bar q\over p}.$$
Thus, $\F(T)^{-1}\xi-x\in L^{\bar q}_{\cF_T}(\O;\dbR^n)$. Then, by
\cite[Theorem 3.1]{Lu-Yong-Zhang 2012}, we can find a $\h v(\cd)\in L^{\bar
q}_\dbF(\O;L^1(0,T;\dbR^n))$ such that
$$\F(T)^{-1}\xi-x=\int_0^T\h v(s)ds.$$
Define
$$v(t)=\F(t)\h v(t),\qq t\in[0,T].$$
Since $q<\bar q$, one has
$$\ba{ll}
\ns\ds\dbE\(\int_0^T|v(t)|dt\)^q\les\dbE\(\int_0^T|\F(t)|\,|\h v(t)|dt\)^q\les\dbE\[\(\sup_{t\in[0,T]}|\F(t)|\)\int_0^T|\h v(t)|dt\]^q\\
\ns\ds\les K\[\dbE\(\int_0^T|\h v(t)|dt\)^{\bar q}\]^{q\over\bar
q}=K\|\h v(\cd)\|_{L^{\bar q}_\dbF(\O;L^1(0.T;\dbR^m))}^{q\over\bar
q}.\ea$$
Now, we define
$$u(t)=B(t)^\top[B(t)B(t)^\top]^{-1}\big[v(t)-A(t)X(t)\big],\qq t\ges0,$$
with $X(\cd)$ defined by \eqref{X(t)}. Then
$$A(t)X(t)+B(t)u(t)=v(t),\qq t\ges0,$$
which implies that
$$X(t)=\F(t)x+\F(t)\int_0^t\F(s)^{-1}\[A(s)X(s)+B(s)u(s)\]ds,\qq t\ges0.$$
This means that $X(\cd)$ is the solution to \eqref{FSDE0},
corresponding to $(x,u(\cd))$ with
$$\ba{ll}
\ns\ds\dbE\(\int_0^T|B(t)u(t)|dt\)^q=K\(\dbE\1n\int_0^T|v(t)-A(t)X(t)|dt\)^q\\
\ns\ds\qq\qq\qq\qq\qq\les
K\[\(\dbE\int_0^T|v(t)|dt\)^q+\(\dbE\int_0^T|X(t)|dt\)^q\].\ea$$
Therefore, $u(\cd)\in\cU^q[0,T]$ which makes $X(T)=\xi$. This proves
our conclusion. \endpf

\ms

Let us make some comments on the above result. To this end, let us
define
$$\dbL_T(u(\cd))=\int_0^Tu(t)dt,\qq u(\cd)\in L^1_\dbF(0,T;\dbR^n).$$
Then a result from \cite[Theorems 3.1, 3.2]{Lu-Yong-Zhang 2012} tells us that
\bel{L=L}\dbL_T\(L^p_\dbF(\O;L^1(0,T;\dbR^n))\)\supseteq\dbL_T\(L^1_\dbF\big(0,T;L^p(\O;\dbR^n)
\big)\)=L^p_{\cF_T}(\O;\dbR^n),\qq\forall p\in[1,\infty),\ee
and
\bel{Lq ne
Lp}\dbL_T\(L^q_\dbF\big(0,T;L^p(\O;\dbR^n)\big)\)\subsetneq
L^p_{\cF_T}(\O;\dbR^n),\qq\forall p\in(1,\infty),~q\in(1,\infty].\ee
Thus,
\bel{Lp ne
Lp}\dbL_T\(L^p_\dbF\big(\O;L^q(0,T;\dbR^n)\big)\)\subsetneq
L^p_{\cF_T}(\O;\dbR^n),\qq\forall 1<p\les q\les\infty.\ee
In particular,
\bel{L2 ne
Lp}\dbL_T\(L^p_\dbF\big(\O;L^2(0,T;\dbR^n)\big)\)\subsetneq
L^p_{\cF_T}(\O;\dbR^n),\qq\forall 1<p\les 2.\ee
Let us look at an implication of the above. Consider a system of the
following form:
$$dX(t)=u(t)dt,\qq t\ges0.$$
For terminal state $\xi\in L^p_{\cF_T}(\O;\dbR^n)$, with $p>1$, if
one is only allowed to use the control from $L^p_\dbF(\O;L^2(0,T;$
$\dbR^n))$, the above system is not controllable. Whereas, by Theorem \ref{D=0}, this system is $L^p$-exactly controllable on $[0,T]$ by $\cU^q[0,T]$ for any $q\in(1,p)$. This is a main reason that we define the $L^p$-exact controllability allowing the control taken from a larger space than
$L^p_\dbF(\O;L^2(0,T;\dbR^m))$, and not restricting the state
process $X(\cd)$ to belong to $L^p_\dbF(\O;C([0,T];\dbR^n))$.

\ms

Next, we notice that in Theorem \ref{D=0}, condition
\eqref{B(t)B(t)>d} implies that
$$\hb{rank}\,B(t)=n\les m,\qq t\in[0,T],~\as$$
This means that the dimension of the control process is no less than
that of the state process. Now, if
\bel{rank B<n}\hb{rank}\,B(t)<n,\qq t\in[0,T],~\as,\ee
which will always be the case if $m<n$, then for each $t\in[0,T]$, there exists a
$\th(t)\in\dbR^n\setminus\{0\}$ such that
\bel{th(t)B(t)=0}\th(t)^\top B(t)=0.\ee
The following gives a negative result for the exact controllability,
under condition \eqref{th(t)B(t)=0} with a little more regularity
conditions on $\th(\cd)$ and $C(\cd)$, which is essentially an
extension of \eqref{Lp ne Lp}.

\bt{non-controllable} \sl Let $D(\cdot)=0$ and {\rm(H1)} hold.
Suppose there exists a continuous differentiable function
$\th:[0,T]\to\dbR^n$, $|\th(t)|=1$, for all $t\in[0,T]$ such that
\eqref{th(t)B(t)=0} holds. Also, let
\bel{C_k}C_k(\cd)\in L^\infty_\dbF(\O;C([0,T];\dbR^{n\times n})),\qq
1\les k\les d.\ee
Then for any $p>1$, $[A(\cd),C(\cd);B(\cd),0]$ is not $L^p$-exactly
controllable on $[0,T]$ by $\cU^p[0,T]$.

\et

\it Proof. \rm Let
$$0=t_0<t_1<t_2<\cds,\qq t_\ell\to T.$$
Define
$$G_0=\bigcup_{\ell=0}^\infty\[t_\ell,{t_\ell+t_{\ell+1}\over2}\),\qq G_1=\bigcup_{\ell=0}^\infty\[{t_\ell+t_{\ell+1}\over2},t_{\ell+1}\).$$
Take
$$\z_0,\z_1\in\dbR^n,\q|\th(T)^\top\z_1-\th(T)^\top\z_0|=1.$$
Set
$$\z(s)=\z_0I_{G_0}(s)-\z_1I_{G_1}(s),\qq s\in[0,T),$$
and
\bel{xi}\xi=x+\int_0^T\sum_{k=1}^d\z(s)dW_k(s).\ee
We claim that the above constructed $\xi$ cannot be hit by the state
$X(T)$ from $X(0)=x$ under any $u(\cd)\in\cU^p[0,T]$. We show this
by contradiction. Suppose otherwise, then for some
$u(\cd)\in\cU^p[0,T]$, $X(\cd)$ satisfies
$$\left\{\2n\ba{ll}
\ns\ds dX(s)=\[A(s)X(s)+B(s)u(s)\]ds+C(s)X(s)dW(s),\qq s\in[0,T],\\
\ns\ds X(0)=x,\qq X(T)=\xi.\ea\right.$$
Hence,
$$d[\th(t)^\top X(t)]=\[\th(t)^\top A(t)+\dot\th(t)^\top\]X(t)dt+\sum_{k=1}^d\th(t)^\top C_k(t)X(t)dW_k(t),\qq t\ges0.$$
Now, let
\bel{eta}\eta(t)=x+\int_0^t\sum_{k=1}^d\z(s)dW_k(s),\qq
t\in[0,T].\ee
Then
\bel{BSDE}\left\{\2n\ba{ll}
\ns\ds d\[\th(t)^\top\(X(t)-\eta(t)\)\]=\[\th(t)^\top A(t)X(t)+\dot\th(t)^\top\(X(t)-\eta(t)\)\]dt\\
\ns\ds\qq\qq\qq\qq\qq\qq+\sum_{k=1}^d\th(t)^\top\(C_k(t)X(t)-\z(t)\)dW_k(t), \quad t\in [0,T],\\
\ns\ds\th(T)^\top\(X(T)-\eta(T)\)=0.\ea\right.\ee
By a standard estimate for BSDEs and Burkholder-Davis-Gundy inequality, one obtains
\bel{SSec3_Eq1}\ba{ll}
\ns\ds\dbE\Big\{\sup_{s\in[t,T]}\Big|\th(s)^\top\(X(s)-\eta(s)\)\Big|^p
+\[\int_t^T\sum_{k=1}^d\Big|\th(s)^\top\(C_k(s)X(s)-\z(s)\)\Big|^2ds\]^{p\over2}\Big\}\\
\ns\ds\les K\dbE\(\int_t^T|\th(s)^\top A(s)X(s)+\dot\th(s)^\top X(s)-\dot\th(s)^\top\eta(s)|ds\)^p\\
\ns\ds\les K(T-t)^p\dbE\[\sup_{s\in[t,T]}|X(s)|^p+\sup_{s\in[t,T]}|\eta(s)|^p\]\\
\ns\ds\les
K(T-t)^p\[|x|^p+\dbE\(\int_0^T|B(s)u(s)|ds\)^p+\(\int_0^T|\z(s)|^2ds\)^{p\over2}\].\ea\ee
On the other hand,
\bel{SSec3_Eq2}\ba{ll}
\ns\ds\dbE\[\sup_{s\in[t,T]}|X(s)-\xi|^p\]\les K\max\{ |T-t|^p,
|T-t|^{p\over2}\}\[\,|x|^p+\dbE\(\int_t^T |B(s)u(s)|ds\)^p\]\\
\ns\ds\qq\qq\qq\qq\qq\qq+K\dbE\(\int_t^T |B(s)u(s)|ds\)^p,\ea\ee
and for any $f,g,h\in L^2(t,T;\dbR)$,
$$\|f-g\|^p_{L^2(t,T;\dbR)}\ges 2^{-(p-1)}\|f-h\|^p_{L^2(t,T;\dbR)}-\|h-g\|^p_{L^2(t,T;\dbR)}.$$
Thus, for each $k=1,2,\cds,d$,
$$\ba{ll}
\ns\ds\dbE\[\int_{t_i}^T\Big|\th(s)^\top\(C_k(s)X(s)-\z(s)\)\Big|^2ds\]^{p\over 2} \\
\ns\ds\ges2^{-4(p-1)}\dbE\(\int_{t_i}^T\big|\th(T)^\top\big[C_k(T)\xi-\z(s)\big]\big|^2 ds\)^{p\over 2}-2^{-3(p-1)}\dbE\(\int_{t_i}^T\big|\big[\th(T)^\top-\th(s)^\top\big]\z(s)\big|^2ds\)^{p\over2}\\
\ns\ds\q-2^{-(p-1)}\dbE\(\int_{t_i}^T\big|\big[\th(s)-\th(T)\big]^\top C_k(T)\xi\Big|^2 ds\)^{p\over 2}
-2^{-2(p-1)}\dbE\(\int_{t_i}^T\big|\th(s)^\top\big[C_k(s)-C_k(T)\big]\xi\big|^2 ds\)^{p\over 2}\\
\ns\ds\q-\dbE\(\int_{t_i}^T\big|\th(s)^\top
C_k(s)\big[X(s)-\xi\big]\big|^2 ds\)^{p\over 2}\ea$$
$$\ba{ll}
\ns\ds\ges 2^{-4(p-1)}\({T-t_i\over2}\)^{p\over 2}\dbE\[
|\th(T)^\top\z_0-\th(T)^\top C_k(T)\xi|^2
+|\th(T)^\top\z_1-\th(T)^\top C_k(T)\xi|^2\]^{p\over 2}~~~~~~~~~~~~~~~~~~~~~~~\\
\ns\ds\qq\qq-2^{-3(p-1)}(T-t_i)^{p\over 2}\|\dot\th(\cd)\|_{L^\infty(0,T;\dbR^n)}^p\(\dbE\int_{t_i}^T|\z(s)|^2ds\)^{p\over2}\\
\ns\ds\qq\qq-2^{-2(p-1)}(T-t_i)^{\frac{p}{2}}\dbE\[\sup_{s\in[t_i,T]}\big|[\th(s)-\th(T)]^\top
C_k(T)\xi\big|^p\]\\
\ns\ds\qq\qq-2^{-(p-1)}(T-t_i)^{p\over 2}\dbE\[\sup_{s\in [t_i,T]}
|C_k(T)-C_k(s)|^p |\xi|^p\]\\
\ns\ds\qq\qq-(T-t_i)^{p\over 2}\|\th(\cd)^\top
C_k(\cd)\|^p_{L^\infty_\dbF(0,T;\dbR^{n\times n})}\dbE\[\sup_{s\in
[t_i,T]}|X(s)-\xi|^p\].\ea$$
Since $C_k(\cd)\in L^\infty_\dbF(\O;C([0,T];\dbR^{n\times n}))$, by
Lebesgue's dominated convergence theorem, we have
$$\dbE\[\sup_{s\in [t_i,T]}|C_k(T)-C_k(s)|^p |\xi|^p
\] = o(1),\qq\as,~~i\to\infty.$$
From \eqref{SSec3_Eq2},
$$\dbE\[\sup_{s\in [t_i,T]}|X(s)-\xi|^p\]=o(1),\qq\as,~~i\to\infty.$$
The other two negative terms can be estimated similarly.
Consequently, making use of \eqref{SSec3_Eq1},
$$\ba{ll}
\ns\ds2^{-4(p-1)}\({T-t_i\over2}\)^{p\over 2}\dbE\[
|\th(T)^\top\z_0-\th(T)^\top C_k(T)\xi|^2
+|\th(T)^\top\z_1-\th(T)^\top C_k(T)\xi|^2 \]^{p\over 2}\\
\ns\ds\les\dbE\(\int_{t_i}^T |\th(s)^\top C_k(s)X(s)
-\th(s)^\top\z(s)|^2 ds \)^{p\over 2} +o\( |T-t_i|^{p\over
2}\)=o\(|T-t_i|^{p\over 2}\) \q \as,~~i\to\infty.\ea$$
This leads to
$$\dbE\(|\th(T)^\top\z_0-\th(T)^\top C_k(T)\xi|^2+|\th(T)^\top\z_1-\th(T)^\top C_k(T)\xi|^2\)^{p\over2}=0.$$
Thus,
$$|\th(T)^\top\z_0-\th(T)^\top C_k(T)\xi|^2+|\th(T)^\top\z_1-\th(T)^\top C_k(T)\xi|^2=0,\qq\as,$$
which is a contradiction since $|\th(T)^\top(\z_0-\z_1)|=1$. The
above implies that the terminal state $\xi$ constructed above cannot
be hit by the state under any $u(\cd)\in\cU^p[0,T]$. This completes
the proof. \endpf

\subsection{The case $D(\cd)$ is surjective}

In this subsection, we let $d=1$. The case $d>1$ can be discussed
similarly. For system $[A(\cd),C(\cd);B(\cd),D(\cd)]$, we assume the
following:
\bel{DD>d}D(t)D(t)^\top\ges\d I,\qq\as,\ae\ t\in[0,T].\ee
In this case, $[D(t)D(t)^\top]^{-1}$ exists and uniformly bounded.
We define
\bel{hABD}\left\{\2n\ba{ll}
\ns\ds\h A(t)=A(t)-B(t)D(t)^\top[D(t)D(t)^\top]^{-1}C(t),\\
\ns\ds\h B(t)=B(t)\big\{I-D(t)^\top[D(t)D(t)^\top]^{-1}D(t)\big\},\\
\ns\ds\h D(t)=B(t)D(t)^\top[D(t)D(t)^\top]^{-1},\ea\right.\ee
and introduce the following controlled system:
\bel{Yu_Sec3_Equation1000}\left\{\2n\ba{ll}
\ns\ds dX(t)=\[\h A(t)X(t)+\h B(t)v(t)+\h D(t)Z(t)\]dt+Z(t)dW(t),\qq t\in [0,T],\\
\ns\ds X(0)=x,\ea\right.\ee
with $X(\cd)$ being the state and $(v(\cd),Z(\cd))$ being the
control. Using our notation, the above system can be denoted by $[\h
A(\cd),0;(\h B(\cd),\h D(\cd)),(0,I)]$. Comparing
$[A(\cd),C(\cd);B(\cd),D(\cd)]$ with $[\h A(\cd),0;(\h B(\cd),\h
D(\cd));(0,I)]$, the latter has a simpler structure. For system $[\h
A(\cd),0;(\h B(\cd),\h D(\cd))$, $(0,I)]$, we need the following
set:
$$\h\cU^p[0,T]\equiv\Big\{v(\cd)\bigm|\h
B(\cd)v(\cd)\in L^p_\dbF(\O;L^1(0,T;\dbR^n))\Big\}.$$
The following result is a kind of reduction.

\bt{Theorem 3.4} \sl Let {\rm (H1)} and \eqref{DD>d} hold. Let $\h
A(\cd),\h B(\cd),\h D(\cd)$ be defined by \eqref{hABD}. Suppose
\bel{|hA|,|hD|}\h A(\cd)\in L^\infty_{\mathbb
F}(\Omega;L^{1+\varepsilon}(0,T;\mathbb R^{n\times n})),\qq \h
D(\cd)\in L^\infty_\dbF(\O;L^2(0,T;\dbR^{n\times n})),\ee
where $\varepsilon>0$ is a given constant. Then system $[\h
A(\cd),0;(\h B(\cd),\h D(\cd)),(0,I)]$ is $L^p$-exactly controllable
on $[0,T]$ by $\h\cU^p[0,T]\1n\times\1n L^p_\dbF(\O;$
$L^2(0,T;\dbR^n))$ if and only if system
$[A(\cd),C(\cd);B(\cd),D(\cd)]$ is $L^p$-exactly controllable on
$[0,T]$ by $\cU^p[0,T]$. \et

\it Proof. \rm ($\Rightarrow$). First of all, we note that for any
$Z(\cd)\in L^p_\dbF(\O;L^2(0,T;\dbR^n))$,
\bel{hDZ}\ba{ll}
\ns\ds\dbE\(\int_0^T|\h D(t)Z(t)|dt\)^p\les\dbE\[\(\int_0^T|\h D(t)|^2dt\)^{p\over2}\(\int_0^T|Z(t)|^2dt\)^{p\over2}\]\\
\ns\ds\les\|\h D(\cd)\|_{L^\infty_\dbF(\O;L^2(0,T;\dbR^{n\times
n}))}^p\dbE\(\int_0^T|Z(t)|^2dt\)^{p\over2}.\ea\ee
Thus, under condition \eqref{|hA|,|hD|}, for any $x\in\dbR^n$ and
$(v(\cd),Z(\cd))\in\h\cU^p[0,T]\times L^p_\dbF(\O;L^2(0,T;\dbR^n))$,
system \eqref{Yu_Sec3_Equation1000} admits a unique solution.

\ms

Now, if system $[\h A(\cd),0;(\h B(\cd),\h D(\cd)),(0,I)]$ is
$L^p$-exactly controllable on $[0,T]$ by $\h\cU^p[0,T]\times
L^p_\dbF(\O;L^2(0,T;$ $\dbR^n))$, then for any $x\in\dbR^n$ and
$\xi\in L^p_{\cF_T}(\O;\dbR^n)$, there exists a triple
$(X(\cd),v(\cd),Z(\cd))\in L^p_\dbF(\O;C([0,T];$
$\dbR^n))\times\h\cU^p[0,T]\times L^p_\dbF(\O;L^2(0,T;\dbR^n))$ such
that
\bel{Yu_Sec3_Equation1}\left\{\2n\ba{ll}
\ns\ds dX(t)=\[\h A(t)X(t)+\h B(t)v(t)+\h
D(t)Z(t)\]dt+Z(t)dW(t),\qq t\in [0,T],\\
\ns\ds X(0)=x,\qq X(T)=\xi.\ea\right.\ee
Define
$$u(t)=D(t)^\top[D(t)D(t)^\top]^{-1}[Z(t)-C(t)X(t)]+[I-D(t)^\top
[D(t)D(t)^\top]^{-1}D(t)]v(t),\quad t\in [0,T].$$
We have
\bel{Yu_Sec3_Bu}B(t)u(t)=\h D(t)[Z(t)-C(t)X(t)]+\h B(t)v(t),\quad
t\in [0,T].\ee
Since $v(\cd)\in\h\cU^p[0,T]$ and
$$\ba{ll}
\ns\ds\dbE\(\int_0^T|\h
D(t)(Z(t)-C(t)X(t))|dt\)^p\les\dbE\[\(\int_0^T|\h D(t)|^2dt
\)^{p\over 2}\(\int_0^T|Z(t)-C(t)X(t)|^2dt\)^{p\over 2}\]\\
\ns\ds\les\|\h D(\cd)\|^p_{L^\infty_\dbF(\O;L^2(0,T;\dbR^{n\times
n}))}
\|Z(\cd)-C(\cd)X(\cd)\|^p_{L^p_\dbF(\O;L^2(0,T;\dbR^n))}<\infty,\ea$$
one has $B(\cd)u(\cd)\in L^p_\dbF(\O;L^1(0,T;\dbR^n))$. Further, by
\bel{Yu_Sec3_Du}D(t)u(t)=Z(t)-C(t)X(t),\quad t\in [0,T],\ee
we obtain $D(\cd)u(\cd)\in L^p_\dbF(\O;L^2(0,T;\dbR^n))$. Therefore
$u(\cd)\in\cU^p[0,T]$. From \eqref{Yu_Sec3_Bu} and
\eqref{Yu_Sec3_Du}, it is easy to see
$$\left\{\2n\ba{ll}
\ns\ds A(t)X(t)+B(t)u(t)=\h A(t)X(t)+\h B(t)v(t)+\h D(t)Z(t),\\
\ns\ds C(t)X(t)+D(t)u(t)=Z(t),\ea\right. \quad t\in [0,T],$$
and thus, \eqref{Yu_Sec3_Equation1} reads
$$\left\{\2n\ba{ll}
\ns\ds dX(t)=\[A(t)X(t)+B(t)u(t)\]dt+\[C(t)X(t)+D(t)u(t)\]dW(t),\qq t\in[0,T],\\
\ns\ds X(0)=x,\qq X(T)=\xi.\ea\right.$$
This proves the $L^p$-exact controllability of system
$[A(\cd),C(\cd);B(\cd),D(\cd)]$ on $[0,T]$ by $\cU^p[0,T]$.

\ms

($\Leftarrow$). If system $[A(\cd),C(\cd);B(\cd),D(\cd)]$ is
$L^p$-exactly controllable on $[0,T]$ by $\cU^p[0,T]$, then for any
$x\in\dbR^n$ and $\xi\in L^p_{\cF_T}(\O;\dbR^n)$, there exists a
pair $(X(\cd),u(\cd))\in
L^p_\dbF(\O;C([0,T];\dbR^n))\times\cU^p[0,T]$ such that
\bel{Yu_Sec3_Equation2}\left\{\2n\ba{ll}
\ns\ds dX(t)=\[A(t)X(t)+B(t)u(t)\]dt+\[C(t)X(t)+D(t)u(t)\]dW(t),\qq t\in[0,T],\\
\ns\ds X(0)=x,\qq X(T)=\xi.\ea\right.\ee
Let
$$
\left\{
\begin{aligned}
& Z(t)=C(t)X(t)+D(t)u(t),\\
& v(t)=u(t),
\end{aligned}
\right. \qq t\in[0,T].
$$
Then $Z(\cd)\in L^p_\dbF(\O;L^2(0,T;\dbR^n))$ and
$$u(t)=D(t)^\top [D(t)D(t)^\top]^{-1}[Z(t)-C(t)X(t)]+[I-D(t)^\top
[D(t)D(t)^\top]^{-1}D(t)]v(t),\quad t\in [0,T].$$
Further,
$$\ba{ll}
\ns\ds B(t)u(t)\1n=\1n B(t)D(t)^\top
[D(t)D(t)^\top]^{-1}[Z(t)\1n-\1n C(t)X(t)]\1n+\1n B(t)[I\1n-\1n
D(t)^\top
[D(t)D(t)^\top]^{-1}D(t)]v(t)\\
\ns\ds\qq\qq=\h D(t)Z(t)+\[\h A(t)-A(t)\]X(t)+\h B(t)v(t),\quad t\in
[0,T].\ea$$
Consequently,
$$\ba{ll}
\ns\ds\dbE\(\int_0^T|\h B(t)v(t)|dt\)^p\les3^{p-1}\dbE\[\(\int_0^T|\h D(t)Z(t)|dt\)^p+\(\int_0^T|[\h A(t)-A(t)]X(t)|dt\)^p\\
\ns\ds\qq\qq\qq\qq\qq\qq+\(\int_0^T|B(t)u(t)|dt\)^p\]\\
\ns\ds\qq\les K\dbE\[\(\int_0^T|\h
D(t)|^2dt\)^{p\over2}\(\int_0^T|Z(t)|^2dt\)^{p\over2}+\sup_{t\in[0,T]}|X(t)|^p
+\(\int_0^T|B(t)u(t)|dt\)^p\]\\
\ns\ds\qq\les
K\[|x|^p+\dbE\(\int_0^T|Z(t)|^2dt\)^{p\over2}+\dbE\(\int_0^T|B(t)u(t)|dt\)^p+\dbE\(\int_0^T|D(t)u(t)|^2dt\)^{p\over2}\].\ea$$
Thus, $v(\cd)\in\h\cU^p[0,T]$. Also,
$$\ba{ll}
\ns\ds\h A(t)X(t)+\h B(t)v(t)+\h D(t)Z(t)\\
\ns\ds=\[A(t)-B(t)D(t)^\top\big[D(t)D(t)^\top\big]^{-1}C(t)\]X(t)
+B(t)\big\{I-D(t)^\top\big[D(t)D(t)^\top\big]^{-1}D(t)\big\}u(t)\\
\ns\ds\qq+B(t)D(t)^\top\big[D(t)D(t)^\top\big]^{-1}\big[C(t)X(t)+D(t)u(t)\big]\\
\ns\ds=A(t)X(t)+B(t)u(t),\quad t\in [0,T].\ea$$
Hence, \eqref{Yu_Sec3_Equation2} can be rewritten as
$$\left\{\2n\ba{ll}
\ns\ds dX(t)=\[\h A(t)X(t)+\h B(t)v(t)+\h
D(t)Z(t)\]dt+Z(t)dW(t),\qq t\in [0,T],\\
\ns\ds X(0)=x,\qq X(T)=\xi.\ea\right.$$
This proves the $L^p$-exact controllability of system $[\h
A(\cd),0;(\h B(\cd),\h D(\cd)),(0,I)]$ on $[0,T]$ by
$\h\cU^p[0,T]\times L^p_\dbF(\O;$ $L^2(0,T;\dbR^n))$. \endpf

\ms

Now, for system $[\h A(\cd),0;(\h B(\cd),\h D(\cd)),(0,I)]$, similar
to the definition of exact controllability, we introduce the
following definition.

\bde{null-controllability} \rm System $[\h A(\cd),0,(\h B(\cd),\h
D(\cd)),(0,I)]$ is said to be {\it exactly null-controllable} by
$\h\cU^p[0,T]\times L^p_\dbF(\O;L^2(0,T;$ $\dbR^n))$ on the time
interval $[0,T]$, if for any $x\in\dbR^n$, there exists a pair
$(v(\cd),Z(\cd))\in\h\cU^p[0,T]\times L^p_\dbF(\O;L^2(0,T;\dbR^n))$
such that the solution $X(\cd)$ to \eqref{Yu_Sec3_Equation1} under
$(v(\cd),Z(\cd))$ satisfies $X(T)=0$. \ede

We have the following result for system $[\h A(\cd),0;(\h B(\cd),\h
D(\cd)),(0,I)]$.

\bt{Liu-Peng} \sl Let {\rm(H1)} and \eqref{DD>d} hold. Let $\h
A(\cd),\h B(\cd),\h D(\cd)$ be defined by \eqref{hABD}. Suppose
\begin{equation}\label{Sec2_ABD_Epsilon}
\begin{aligned}
& \widehat A(\cdot) \in L^\infty_{\mathbb
F}(\Omega;L^{1+\varepsilon}(0,T;\mathbb R^{n\times n})),\quad
\widehat B(\cdot) \in L^{(2\vee p) +\varepsilon}_{\mathbb
F}(\Omega;L^2(0,T;\mathbb R^{n\times m})),\\
& \widehat D(\cdot)\in L^\infty_{\mathbb
F}(\Omega;L^{2+\varepsilon}(0,T;\mathbb R^{n\times n})),
\end{aligned}
\end{equation}
where $2\vee p \equiv \max\{2,p\}$ and $\varepsilon>0$ is a given
constant. Then the following are equivalent:\\

\ss

{\rm(i)} $[\h A(\cd),0;(\h B(\cd),\h D(\cd)),(0,I)]$ is
$L^p$-exactly controllable on $[0,T]$ by $\h\cU^p[0,T]\times
L^p_\dbF(\O;L^2(0,T;\dbR^n))$;\\

\ss

{\rm(ii)} $[\h A(\cd),0;(\h B(\cd),\h D(\cd)),(0,I)]$ is exactly
null-controllable on $[0,T]$ by $\h\cU^p[0,T]\times
L^p_\dbF(\O;L^2(0,T;\dbR^n))$;\\

\ss

{\rm(iii)} Matrix $G$ defined below is invertible:
\bel{Yu_Sec3_G}G=\dbE\int_0^T\cY(t)\h B(t)\h B(t)^\top \cY(t)^\top
dt,\ee
where $\cY(\cd)$ is the adapted solution to the following FSDE:
\bel{cY}\left\{\2n\ba{ll}
\ns\ds d\cY(t)=-\cY(t)\h A(t) dt-\cY(t)\h D(t)dW(t),\qq t\ges0,\\
\ns\ds\cY(0)=I.\ea\right.\ee

\et

\it Proof. \rm (i) $\Ra$ (ii) is trivial.

\ms

(ii) $\Ra$ (iii). First of all, under \eqref{Sec2_ABD_Epsilon}, FSDE
\eqref{cY} admits a unique strong solution $\cY(\cdot)\in
\bigcap_{q>1}L^q_{\mathbb F}(\Omega;C([0,T];$ $\dbR^{n\times n}))$
and the Matrix $G$ is well defined. Now we prove the conclusion by
contradiction. Suppose matrix $G$ is not invertible, then there
exists a vector $0\ne\b\in\dbR^n$ such that
$$0=\b^\top G\b=\dbE\int_0^T\b^\top\cY(t)\h B(t)\h B(t)^\top\cY(t)^\top\b dt=\dbE\int_0^T|\h B(t)^\top\cY(t)^\top\b|^2dt.$$
Therefore
\bel{bYB=0}\b^\top\cY(t)\h B(t) =0,\qq\ae\;t\in[0,T],\q\as\ee
Now, we claim that by choosing $x=\b\in\dbR^n$, there will be no
$(v(\cd),Z(\cd))\in\h\cU^p[0,T]\times L^p_\dbF(\O;L^2(0,T;$
$\dbR^n))$ such that the corresponding state precess $X(\cd)\equiv
X(\cd\,;x,v(\cd),Z(\cd))$ satisfies
$$X(0)=x,\qq X(T)=0$$
In fact, suppose there exists a pair
$(v(\cd),Z(\cd))\in\h\cU^p[0,T]\times L^p_\dbF(\O;L^2(0,T;\dbR^n))$
such that the above is true. Then applying the It\^o's formula to
$\cY(t)X(t)$ on the interval $[0,T]$, one obtains the following
relationship:
\bel{Yu_Sec3_ItoFormulaHat}-\b=\cY(T)X(T)-X(0)=\int_0^T\cY(t)\h
B(t)v(t)dt+\int_0^T\[\cY(t)Z(t) -\cY(t)\h D(t)X(t)\]dW(t).\ee
It is easy to check that
$$\dbE\(\int_0^T|\cY(t)Z(t)-\cY(t)\h D(t)X(t)|^2dt\)^{1/2}<\infty.$$
Thus,
\bel{Yu_Sec3_xHat}-\b=\dbE\int_0^T\cY(t)\h B(t)v(t)dt.\ee
Making use of \eqref{bYB=0}, we get
$$-|\b|^2=\dbE\int_0^T\b^\top\cY(t)\h B(t)v(t)dt =0,$$
a contradiction.

\ms

(iii) $\Ra$ (i). Under \eqref{Sec2_ABD_Epsilon}, for any given
$\xi\in L^p_{\cF_T}(\O;\dbR^n)$, the following BSDE
$$\left\{\2n\ba{ll}
\ns\ds dX_1(t)=\[\h A(t)X_1(t) +\h D(t)Z_1(t)\]dt +Z_1(t)dW(t),\quad t\in [0,T],\\
\ns\ds X_1(T)=\xi\ea\right.$$
admits a unique adapted solution $(X_1(\cd),Z_1(\cd))\in
L^p_\dbF(\O;C([0,T];\dbR^n))\times L^p_\dbF(\O;L^2(0,T;\dbR^n))$.
Since $G$ is invertible, for any $x\in\dbR^n$, we may define
$$v(t)=-\h B(t)^\top\cY(t)^\top G^{-1}\big[x-X_1(0)\big],\quad t\in [0,T].$$
Note that $\h B(\cd)\in L^{(2\vee
p)+\varepsilon}_\dbF(\O;L^2(0,T;\dbR^{n\times m}))$ leads to the
following:
$$\ba{ll}
\ns\ds\dbE\(\int_0^T|\h B(t)v(t)|dt\)^p\les K\dbE\(\int_0^T|\h B(t)|^2|\cY(t)|dt\)^p\les K\dbE\[\(\sup_{t\in[0,T]}|\cY(t)|^p\)\(\int_0^T|\h B(t)|^2dt\)^p\]\\
\ns\ds\les K\[\dbE\(\int_0^T|\h
B(t)|^2dt\)^{p+\varepsilon}\]^{\frac{p}{p+\varepsilon}}\[\dbE\(\sup_{t\in
[0,T]}|\cY(t)|^{\frac{p\varepsilon}{
p+\varepsilon}}\)\]^{\frac{\varepsilon }{
p+\varepsilon}}<\infty,\ea$$
which implies $v(\cd)\in\h\cU^p[0,T]$. For this $v(\cd)$, we define
$(X_2(\cd),Z_2(\cd))\in L^p_\dbF(\O;C([0,T];\dbR^n))\times
L^p_\dbF(\O;L^2(0,T;\dbR^n))$ to be the unique adapted solution of
the following BSDE:
$$\left\{\2n\ba{ll}
\ns\ds dX_2(t)=\[\h A(t)X_2(t)+\h B(t)v(t)+\h D(t)Z_2(t)\]dt+Z_2(t)dW(t),\q t\in [0,T],\\
\ns\ds X_2(T)=0.\ea\right.$$
Applying It\^o's formula to $\cY(\cd)X_2(\cd)$, we have (comparing
with \eqref{Yu_Sec3_ItoFormulaHat})
$$\ba{ll}
\ns\ds-X_2(0)=\cY(T)X_2(T)-X_2(0)\\
\ns\ds=\int_0^T\cY(t)\h B(t)v(t)dt+\int_0^T\[\cY(t)Z_2(t)-\cY(t)\h D(t)X_2(t)\]dW(t)=\dbE\int_0^T\cY(t)\h B(t)v(t)dt\\
\ns\ds=-\Big\{\int_0^T\dbE\[\cY(t)\h B(t)\h B(t)^\top
\cY(t)^\top\]dt\Big\}G^{-1}\big[x-X_1(0)\big]=-\big[x-X_1(0)\big],\ea$$
which implies $X_1(0)+X_2(0)=x$. Now, we define
$$X(t)=X_1(t)+X_2(t),\qq Z(t)=Z_1(t)+Z_2(t),\qq t\in [0,T].$$
Then, by linearity, we see that $(X(\cd),v(\cd),Z(\cd))$ satisfies
the following:
$$\left\{\2n\ba{ll}
\ns\ds dX(t)=\[\h A(t)X(t)+\h B(t)v(t)+\h D(t)Z(t)\]dt+Z(t)dW(t),\q t\in [0,T],\\
\ns\ds X(0)=x,\qq X(T)=\xi.\ea\right.$$
This means that system $[\h A(\cd),0;(\h B(\cd),\h D(\cd)),(0,I)]$
is $L^p$-exactly controllable on $[0,T]$ by $\h\cU^p[0,T]\times
L^p_\dbF(\O;L^2$ $(0,T;\dbR^n))$. \endpf

\ms

The above result is essentially due to Liu--Peng \cite{Liu-Peng
2010}. We have re-organized the way presenting the result. It is
worthy of pointing out that, unlike in \cite{Liu-Peng 2010}, we
allow the coefficients to be unbounded and allow $p$ to be different
from $2$. Combining Theorems \ref{Theorem 3.4} and \ref{Liu-Peng},
we have the following result.

\bt{Theorem 3.7} \sl Let {\rm (H1)}, \eqref{DD>d} and
\eqref{Sec2_ABD_Epsilon} hold with $\h A(\cd)$, $\widehat B(\cdot)$
and $\h D(\cd)$ defined by \eqref{hABD}. Then system
$[A(\cd),C(\cd);B(\cd),D(\cd)]$ is $L^p$-exactly controllable on
$[0,T]$ by $\cU^p[0,T]$ if and only if $G$ defined by
\eqref{Yu_Sec3_G} is invertible.
\end{theorem}

As a simple corollary of the above, we have the following result for
the case of deterministic coefficients.

\bc{deterministic coef} \sl Let {\rm(H1)}, \eqref{DD>d} and
\eqref{Sec2_ABD_Epsilon} hold. Let $\h A(\cd)$ and $\h B(\cd)$ be
deterministic. Let $\h\F(\cd)$ be the solution to the following ODE:
\bel{h Phi}\left\{\2n\ba{ll}
\ns\ds d\h\F(t)=\h A(t)\h\F(t)dt,\qq t\ges0,\\
\ns\ds\h\F(0)=I.\ea\right.\ee
Denote
$$\Psi=\int_0^T\h\Phi(s)^{-1}\h B(s)\(\h\Phi(s)^{-1}\h B(s)\)^\top ds.$$
Suppose that $\Psi$ is invertible. Then system
$[A(\cd),C(\cd);B(\cd),D(\cd)]$ is $L^p$-exactly controllable on
$[0,T]$ by $\cU^p[0,T]$. \ec

\it Proof. \rm In the current case, we have
$$\h\F(t)^{-1}=\dbE\cY(t),\qq t\in[0,T].$$
On the other hand,
$$\ba{ll}
\ns\ds0\les\dbE\[\(\cY(t)-\dbE\cY(t)\)\h B(t)\h B(t)^\top\(\cY(t)-\dbE\cY(t)\)^\top\]\\
\ns\ds\q=\dbE\[\cY(t)\h B(t)\h B(t)^\top\cY(t)\]-\[\dbE\cY(t)\h
B(t)\h B(t)^\top\dbE\cY(t)\].\ea$$
Hence,
$$G\ges\Psi.$$
Then our conclusion follows from the above theorem. \endpf

\ms

The invertibility of matrix $G$ defined by \eqref{Yu_Sec3_G} gives a
nice criterion for the $L^p$-exact controllability of system
$[A(\cd),C(\cd);$ $B(\cd),D(\cd)]$ through the $L^p$-exact
controllability of system $[\h A(\cd),0;(\h B(\cd),\h
D(\cd)),(0,I)]$. However, unless $n=1$, in the case of random
coefficients, the solution $\cY(\cd)$ of FSDE \eqref{cY} does not
have a relatively simple (explicit) form. Thus, the applicability of
condition (iii) in Theorem \ref{Liu-Peng} is somehow limited. In the
rest of this subsection, we will present another sufficient
condition for the $L^p$-exact controllability of $[\h A(\cd),0;(\h
B(\cd),\h D(\cd)),(0,I)]$ which might have a better applicability.

\ms

Now, with random coefficients, we still let $\h\F(\cd)$ be the
solution to \eqref{h Phi} which is a random ODE. Presumably,
$\h\F(\cd)$ is easier to get than $\cY(\cd)$ (the solution of
\eqref{cY}). Define
$$\wt D(t)=\h\F(t)^{-1}\h D(t)\h\F(t),\qq t\in[0,T],$$
and introduce the following {\it mean-field stochastic
Fredholm integral equation} of first kind:
\bel{MF-Fredholm}Y(t)=\z+\L\int_\t^{\bar\t}\dbE_\t\big[\wt D(s)\h Z(s)\big]ds -\int_t^{\bar\t}\wt D(s)\h Z(s)ds-\int_t^{\bar\t}\h Z(s)dW(s),\qq
t\in[\t,\bar\t],\ee
where $\dbE_\t[\,\cd\,]=\dbE[\,\cd\,|\cF_\t]$. We have the following result.

\bl{Yu_Sec3_Lemma_MF_L} \sl Let the following hold:
\bel{|hA|,|hD|g}\h A(\cd)\in L^\infty_\dbF(\Omega;L^{1+\varepsilon}(0,T;\dbR^{n\times n})),\qq
\h D(\cd)\in L^\infty_\dbF(\O;L^{2+\varepsilon}(0,T;\dbR^{n\times n})),\ee
for some $\varepsilon>0$. Then  there exists a positive constant
$\e'$ depending only on $\h D(\cdot)$ and $\L$, such that for any
$0\les\t<\bar\t\les T$, any $\z\in L^p_{\cF_{\bar\t}}(\O;\dbR^n)$
and $\L\in L^\infty_{\cF_{\bar\t}}(\O;\dbR^{n\times n})$, satisfying
$\bar\tau-\tau\les \e'$, \eqref{MF-Fredholm} admits a unique solution
$(Y(\cd),\h Z(\cd))\in L^p_\dbF(\O;C([\t,\bar\t];\dbR^n))\times L^p_\dbF(\O;L^2(\t,\bar\t;\dbR^n))$.

\el

\it Proof. \rm  Let $0\les\t<\bar\t\les T$, and let $\h z(\cd)\in L^p_\dbF(\O;L^2(\t,\bar\t;\dbR^n))$ be given and consider the following BSDE:
$$Y(t)=\z+\L\int_\t^{\bar\t}\dbE_\t\big[\wt D(s)\h z(s)\big]ds-\int_t^{\bar\t}\wt
D(s)\h Z(s)ds-\int_t^{\bar\t}\h Z(s)dW(s),\qq t\in[\t,\bar\t].$$
On the right hand side of the above, the sum of the first two terms
are treated as the terminal state. By the standard theory of BSDEs
(\cite{El Karoui-Peng-Quenez 1997}), we know that the above BSDE
admits a unique adapted solution $(Y(\cd),\h Z(\cd))$ and
the following estimate holds:
$$\dbE\[\sup_{t\in[\t,\bar\t]}|Y(t)|^p+\(\int_\t^{\bar\t}|\h Z(s)|^2ds\)^{p\over2}\]\les K\dbE\Big|\z+\L\int_\t^{\bar\t}\dbE_\t\big[\wt D(s)
\h z(s)\big]ds\Big|^p.$$
Thus, for $\h z_1(\cd),\h z_2(\cd)\in
L^p_\dbF(\O;L^2(\t,\bar\t;\dbR^n))$, if we let $(Y_1(\cd),\h
Z_1(\cd))$ and $(Y_2(\cd),\h Z_2(\cd))$ be the corresponding adapted
solutions, then, noting that both $\h\F(\cd)$ and $\h\F(\cd)^{-1}$
are bounded, and $\h D(\cd)\in
L^\infty_\dbF(\O;L^{2+\varepsilon}(0,T;\dbR^{n\times n}))$, one has
$$\ba{ll}
\ns\ds\dbE\[\sup_{t\in [\t,\bar\t]}|Y_1(t)-Y_2(t)|^p
+\(\int_\t^{\bar\t}|\h Z_1(s)-\h Z_2(s)|^2ds\)^{p\over2}\]\\
\ns\ds\les K\dbE\Big|\int_\t^{\bar\t}\dbE_\tau\[|\h D(s)|\,|\h z_1(s)-\h z_2(s)|\]ds \Big|^p\les K\dbE\Big\{\int_\t^{\bar\t}\[|\h D(s)|\,|\h z_1(s)-\h z_2(s)|\]ds \Big\}^p\\
\ns\ds\les K\dbE\Big\{\(\int_\t^{\bar\t}|\h D(s)|^2ds\)^{1\over2}\( \int_\t^{\bar\t}|\h z_1(s)-\h z_2(s)|^2ds\)^{1\over2}\Big\}^p\\
\ns\ds\les K \sup_{\omega\in \Omega}\(\int_\t^{\bar\t}|\h D(s,\omega)|^2ds\)^{p\over2}
\dbE\(\int_\t^{\bar\t}|\h z_1(s)-\h z_2(s)|^2ds\)^{p\over2}\\
\ns\ds\les K(\bar\t-\t)^{\e
p\over2(2+\e)}\[\sup_{\omega\in\Omega}\(\int_0^{T}|\h
D(s,\omega)|^{2+\e} ds\)^{p\over2+\e}\]\dbE\(\int_\t^{\bar\t}|\h
z_1(s)-\h z_2(s)|^2ds \)^{p\over2}.\ea$$
Consequently, we can find an absolute constant $\varepsilon'>0$ such
that as long as $0<\bar\t-\t\les\varepsilon'$, the map $\h
z(\cd)\mapsto\h Z(\cd)$ is a contraction which admits a unique fixed
point $\h Z(\cd)$. Letting $Y(\cd)$ be given by \eqref{MF-Fredholm},
one sees that $(Y(\cd),\h Z(\cd))\in
L^p_\dbF(\O;C([\t,\bar\t];\dbR^n))\times L^p_\dbF(\t,\bar\t;\dbR^n)$
is the unique solution to \eqref{MF-Fredholm}. That completes the
proof.
\endpf

\ms

Now, we are ready to prove the following result.

\bt{Sufficient 1} \sl Let {\rm(H1)},
\eqref{DD>d} and \eqref{Sec2_ABD_Epsilon} hold. Let
\bel{Psi,Th}\left\{\2n\ba{ll}
\ns\ds\Psi(t,\t)=\int_\t^t\[\dbE_\t\(\h\F(s)^{-1}\h B(s)\)\dbE_\t\([\h\F(s)^{-1}\h B(s)]^\top\)\]ds,\\
\ns\ds\Th(t,\t)=\int_\t^t\[\h\F(s)^{-1}\h B(s)\dbE_\t\([\h\F(s)^{-1}\h B(s)]^\top\)\]ds,\ea\right.\qq0\les\t<t\les T.\ee
Suppose there exists a $\d>0$ such that, for any $T-\d\les\t<t\les T$,
$\Psi(t,\t)$ is invertible and $\Psi(t,\t)^{-1} \in
L^\infty_{\cF_\t}(\O;\dbR^{n\times n})$. Moreover,
suppose there exists a constant $M>0$ such that
\bel{M}|\Th(t,\t)\Psi(t,\t)^{-1}|\les M,\qq T-\d\les\t<t\les
T,\q\as\ee
Then, for any $p>1$, system $[A(\cd),C(\cd);B(\cd),D(\cd)]$ is
$L^p$-exactly controllable on $[0,T]$ by $\cU^p[0,T]$.

\et

\ms

\rm

Note that for $\Psi(t,\t)$ and $\Th(t,\t)$ defined in
\eqref{Psi,Th}, one has
$$\dbE_\t\Th(t,\t)=\Psi(t,\t).$$
Thus, in the case that $\h A(\cd)$ and $\h B(\cd)$ are
deterministic, condition \eqref{M} is automatically true with $M=1$,
as long as $\Psi(t,\t)$ is invertible. We will present an example
that $\h A(\cd)$ is random and \eqref{M} holds. We also point out
that condition \eqref{DD>d} implies that $m\ges n$. Further,
condition that $\Psi(t,\t)^{-1}$ exists implies that $m>n$. In fact,
if \eqref{DD>d} holds and $m=n$, then $\h B(\cd)=0$ which implies
that $\Psi(t,\t)=0$. We will say a little bit more about this
shortly.

\ms

\it Proof. \rm By Theorem \ref{Theorem 3.4}, all we need do is to
prove the $L^p$-exact controllability of system $[\h A(\cd),0;(\h
B(\cd),$ $\h D(\cd)),(0,I)]$ on $[0,T]$ by $\h\cU^p[0,T]\times
L^p_\dbF(\O;L^2(0,T;\dbR^n))$. Now, let $T-\d\les\t<\bar\t\les T$,
and $\xi\in L^p_{\cF_\t}(\O;\dbR^n)$, $\bar\xi\in
L^p_{\cF_{\bar\t}}(\O;\dbR^n)$. For any given $v(\cd),Z(\cd)$, the
solution to $[\h A(\cd),0;(\h B(\cd),\h D(\cd)),(0,I)]$ on
$[\t,\bar\t]$ with $X(\t)=\xi$ is given by
\bel{X2(t)}X(t)=\h\F(t)\h\F(\t)^{-1}\xi+\h\F(t)\int_\t^t\h\F(s)^{-1}\[\h
B(s)v(s) +\h D(s)Z(s)\]ds+\h\F(t)\int_\t^t\h\F(s)^{-1}Z(s)dW(s).\ee
Therefore, getting $X(\bar\t)=\bar\xi$ is equivalent to having the following:
\bel{X(T)=xi}\h\F(\bar\t)^{-1}\bar\xi-\h\F(\t)^{-1}\xi=\int_\t^{\bar\t}\h\F(s)^{-1}\[\h B(s)v(s)+\h D(s)Z(s)\]ds+\int_\t^{\bar\t}\h\F(s)^{-1}Z(s)dW(s).\ee
This implies
\bel{3.22}\dbE_\t[\h\F(\bar\t)^{-1}\bar\xi]-\h\F(\t)^{-1}\xi-\int_\t^{\bar\t}\dbE_\t\[\h\F(s)^{-1}\h D(s)Z(s)\]ds=\int_\t^{\bar\t}\dbE_\t\[\h\F(s)^{-1}\h B(s)v(s)\]ds.\ee
We now take
$$\ba{ll}
\ns\ds v(s)=\dbE_\t\(\h\F(s)^{-1}\h B(s)\)^\top\Psi(\bar\t,\t)^{-1}\[\dbE_\t[\h\F(\bar\t)^{-1}\bar\xi]-\h\F(\t)^{-1}\xi\\
\ns\ds\qq\qq\qq\qq\qq\qq\qq\qq\qq-\int_\t^{\bar\t}\dbE_\t\(\h\F(t)^{-1}\h D(t)Z(t)\)dt\],\qq s\in[\t,\bar\t],\ea$$
which is $\cF_\t$-measurable. Further, noting that $\h\F(\cd)^{-1}$
is bounded and $\Psi(\bar\t,\t)^{-1}\in L^\infty_{\cF_\t}(\O;\dbR^{n\times n})$, one has
$$\dbE\(\int_\t^{\bar\t}|\h B(s)v(s)|ds\)^p\les K\dbE\Big\{\(\int_\t^{\bar\t}|\h B(s)|\dbE_\t|\h B(s)|ds\)^p\G\Big\},$$
where
$$\G\equiv\dbE_\t|\bar\xi|^p+|\xi|^p+\(\int_\t^{\bar\t}\dbE_\t\big[\,|\h D(t)||Z(t)|\,\big]dt\)^p,$$
which is $\cF_\t$-measurable. Since
$$\ba{ll}
\ns\ds\(\int_\t^{\bar\t}|\h B(s)|\dbE_\t|\h B(s)|ds\)^p\les\({1\over2}\int_\t^{\bar\t}
|\h B(s)|^2ds+{1\over2}\int_\t^{\bar\t}\dbE_\t|\h B(s)|^2ds\)^p\\
\ns\ds\qq\qq\qq\qq\qq\q\les K\(\int_\t^{\bar\t}|\h B(s)|^2ds\)^p+K\dbE_\t\(\int_\t^{\bar\t}|\h B(s)|^2ds\)^p,\ea$$
we have
$$\ba{ll}
\ns\ds\dbE\(\int_\t^{\bar\t}|\h B(s)v(s)|ds\)^p\les K\dbE\Big\{\(\int_\t^{\bar\t}
|\h B(s)|^2ds\)^p\G+\dbE_\t\(\int_\t^{\bar\t}|\h B(s)|^2ds\)^p\G\Big\}\\
\ns\ds=K\dbE\Big\{\(\int_\t^{\bar\t}|\h B(s)|^2ds\)^p\G\Big\}+K\dbE\Big\{\dbE_\t\[\( \int_\t^{\bar\t}|\h B(s)|^2ds\)^p\G\]\Big\}\les K\dbE\Big\{\(\int_\t^{\bar\t}|\h B(s)|^2ds\)^p\G\Big\}.\ea$$
Since $\h B(\cd)\in L^\infty_\dbF(\O;L^2(0,T;\dbR^{n\times m}))$ and
$\h D(\cd)\in L^\infty_\dbF(\O;L^2(0,T;\dbR^{n\times m}))$, one has
$$\ba{ll}
\ns\ds\dbE\(\int_\t^{\bar\t}|\h B(s)v(s)|ds\)^p\les K\dbE\Big\{|\xi|^p+|\bar\xi|^p+\[
\int_\t^{\bar\t}\dbE_\t\(|\h D(t)||Z(t)|\)dt\]^p\Big\}\\
\ns\ds\les K\dbE\(|\xi|^p+|\bar\xi|^p\)+K\dbE\Big\{\dbE_\t\[\(\int_\t^{\bar\t}|\h
D(t)|^2dt\)^{1\over2}\(\int_\t^{\bar\t}|Z(t)|^2dt\)^{1\over2}\]\Big\}^p\\
\ns\ds\les K\dbE\(|\xi|^p+|\bar\xi|^p\)+K\dbE\[\dbE_\t\(\int_\t^{\bar\t}|Z(t)|^2dt
\)^{1\over2}\]^p\les K\dbE\[|\xi|^p+|\bar\xi|^p+\(\int_\t^{\bar\t}|Z(t)|^2dt\)^{p\over2}\].\ea$$
Thus, $Z(\cd)\in L^p_\dbF(\O;L^2(\t,\bar\t;\dbR^n))$ implies $\h
B(\cd)v(\cd)\in L^p_\dbF(\O;L^1(\t,\bar\t;\dbR^n))$. For such a
$v(\cd)$, \eqref{3.22} holds. Moreover, \eqref{X(T)=xi} becomes
\bel{Fredholm}0=\eta(\bar\t,\t)+\int_\t^{\bar\t}\Th(\bar\t,\t)\Psi(\bar\t,\t)^{-1}\dbE_\t\(\wt
D(s)\h Z(s)\)ds -\int_\t^{\bar\t}\wt D(s)\h Z(s)ds-\int_\t^{\bar\t}\h Z(s)dW(s),\ee
where
$$\wt D(s)=\h\F(s)^{-1}\h D(s)\h\F(s),\qq\h Z(s)=\h\F(s)^{-1}Z(s),$$
and
$$\eta(\bar\t,\t)\equiv\h\F(\bar\t)^{-1}\bar\xi-\h\F(\t)^{-1}\xi-\Th(\bar\t,\t)
\Psi(\bar\t,\t)^{-1}\(\dbE_\t[\h\F(\bar\t)^{-1}\bar\xi]-\h\F(\t)^{-1}\xi\),$$
which is $\cF_{\bar\t}$-measurable with
$$\dbE|\eta(\bar\t,\t)|^p\les K\dbE\(|\xi|^p+|\bar\xi|^p\),\qq\dbE_\t\eta(\bar\t,\t)=0.$$
To summarize the above, we see that for given $0\les\t<\bar\t\les T$
and $\xi\in L^p_{\cF_\t}(\O;\dbR^n)$, $\bar\xi\in
L^p_{\cF_{\bar\t}}(\O;\dbR^n)$, to get a pair of $(v(\cd),Z(\cdot))$
so that $X(\t)=\xi$ and $X(\bar\t)=\bar\xi$ if and only if there
exists an $\dbF$-adapted process $\h Z(\cd)$ such that
\eqref{Fredholm} holds true.
%

\ms

By Lemma \ref{Yu_Sec3_Lemma_MF_L}, there exists a uniform
$\varepsilon'>0$ such that as long as $ (T-\delta)\vee
(T-\varepsilon')\les\tau<\bar\tau\les T$, the following equation
\bel{BSDE*}\ba{ll}
\ns\ds Y(t)=\eta(\bar\t,\t)+\int_\t^{\bar\t}\Th(\bar\t,\t)\Psi(\bar\t,\t)^{-1}\dbE_\t\(\wt D(s)\h Z(s)\)ds-\int_t^{\bar\t}\wt D(s)\h Z(s)ds\\
\ns\ds\qq\qq\qq\qq\qq\qq\qq\qq -\int_t^{\bar\t}\h Z(s)dW(s),\qq
t\in[\t,\bar\t]\ea\ee
admits a unique adapted solution $(Y(\cd),\h Z(\cd))\in
L^p_\dbF(\O;C([\t,\bar\t];\dbR^n))\times L^p_\dbF(\O;L^2(\t,\bar\t;\dbR^n))$.
We notice that
$$Y(\t)=\dbE_\t[Y(\t)]=0,$$
which proves the existence of a $\h Z(\cd)\in L^p_\dbF(\t,\bar\t;\dbR^n)$ to \eqref{Fredholm}.

\ms

Now, we can complete our proof as follows: Arbitrarily select a $T_0$ such that $(T-\d)\vee(T-\e')\les T_0<T$. For any $x\in\dbR^n$ and $\bar\xi\in L^p_{\cF_T}(\O;\dbR^n)$, let
$$\xi=X(T_0;x,0,0)\in L^p_{\cF_{T_0}}(\O;\dbR^n).$$
Then by what we have proved, there exists a pair
$(v(\cd),Z(\cd))\in\h\cU^p[T_0,T]\times
L^p_\dbF(\O;L^2(T_0,T;\dbR^n))$ such that
$$X(T_0)=\xi,\qq X(T)=\bar\xi.$$
We obtain the $L^p$-exact controllability on $[0,T]$ by
$\cU^p[0,T]$. \endpf

\ms

The following simple example is to show that condition \eqref{M} is possible for random coefficient case.

\bex{} \rm  Let
$$\h A(t)=\begin{pmatrix}0&a(t)\\ 0&0\end{pmatrix},\qq\h B(t)=\begin{pmatrix}0\\ 1\end{pmatrix},$$
where $a:[0,T]\times\O\to\dbR$ satisfies the following conditions:
$t\mapsto a(t)$ is $C^2$ with
$$a(\cd),a'(\cd),a''(\cd)\in L^\infty_\dbF(0,T;\dbR),\qq a(t)\ges1,\qq t\in[0,T],~\as$$
For example, we may choose
$$a(t)=1+\int_0^t\int_0^s{W(\t)^2\over1+W(\t)^2}d\t ds,\qq t\in[0,T].$$
Then
$$\h\F(t)=\begin{pmatrix}1&\int_0^ta(s)ds\\ 0&1\end{pmatrix},\qq t\in[0,T].$$
Thus,
$$\h\F(t)^{-1}\h B(t)=\begin{pmatrix}1&-\int_0^ta(s)ds\\ 0&1\end{pmatrix}\begin{pmatrix}0\\ 1\end{pmatrix}=\begin{pmatrix}-\int_0^ta(s)ds\\ 1\end{pmatrix}\equiv\begin{pmatrix}\a(t)\\ 1\end{pmatrix}.$$
Denoting $\bar\a(s)=\dbE_\t\a(s)$, we have
$$\ba{ll}
\ns\ds\Psi(t,\t)=\int_\t^t\dbE_\t\(\h\F(s)^{-1}\h
B(s)\)\dbE_\t\(\h\F(s)^{-1}\h B(s)\)^\top ds
=\int_\t^t\begin{pmatrix}\bar\a(s)\\ 1\end{pmatrix}\begin{pmatrix}\bar\a(s)&1\end{pmatrix}ds\\
\ns\ds\qq\q=\begin{pmatrix}\int_\t^t\bar\a(s)^2ds&\int_\t^t\bar\a(s)ds\\
\int_\t^t\bar\a(s)ds&t-\t\end{pmatrix}.\ea$$
A direct computation shows that (denoting $\bar a(\t)=\dbE_\t a(\t)$)
$$F(t)\equiv\det\Psi(\t,t)=(t-\t)\int_\t^t\bar\a(s)^2ds-\(\int_\t^t\bar\a(s)ds\)^2=\[2\bar a(\t)^2+R(t)\](t-\t)^4,$$
with
$$\lim_{t\to\t}R(t)=0.$$
Hence, for $t-\t>0$ small, $\Psi(\t,t)$ is invertible, and
$$\Psi(t,\t)^{-1}={1\over F(t)}\begin{pmatrix}t-\t&-\int_\t^t\bar\a(s)ds\\-\int_\t^t\bar\a(s)ds&\int_\t^t\bar\a(s)^2ds
\end{pmatrix}.$$
Also,
$$\ba{ll}
\ns\ds\Th(t,\t)=\int_\t^t\h\F(s)^{-1}\h B(s)\dbE_\t\(\h\F(s)^{-1}\h
B(s)\)^\top ds
=\int_\t^t\begin{pmatrix}\a(s)\\ 1\end{pmatrix}\begin{pmatrix}\bar\a(s)&1\end{pmatrix}ds\\
\ns\ds\qq\q=\begin{pmatrix}\int_\t^t\a(s)\bar\a(s)ds&\int_\t^t\a(s)ds\\
\int_\t^t\bar\a(s)ds&t-\t\end{pmatrix}. \ea$$
Then
$$\ba{ll}
\ns\ds\Th(t,\t)\Psi(t,\t)^{-1}={1\over F(t)}\begin{pmatrix}\int_\t^t\a(s)\bar\a(s)ds&\int_\t^t\a(s)ds\\
\int_\t^t\bar\a(s)ds&t-\t\end{pmatrix}\begin{pmatrix}t-\t&-\int_\t^t\bar\a(s)ds\\-\int_\t^t\bar\a(s)ds&\int_\t^t\bar\a(s)^2ds
\end{pmatrix}\\
\ns\ds\qq\qq\qq\q\equiv{1\over
F(t)}\begin{pmatrix}\L_1(t)&\L_2(t)\\
0&F(t)\end{pmatrix},\ea$$
where
$$\ba{ll}
\ns\ds\L_1(t)=(t-\t)\int_\t^t\a(s)\bar\a(s)ds-\(\int_\t^t\a(s)ds\)\(\int_\t^t\bar\a(s)ds\),\\
\ns\ds\L_2(t)=\(\int_\t^t\a(s)ds\)\(\int_\t^t\bar\a(s)^2ds\)-\(\int_\t^t\bar\a(s)ds\)\(\int_\t^t\a(s)\bar\a(s)ds\).
\ea$$
Some direct (lengthy) calculations show that
$$\ba{ll}
\ns\ds\L_1(t)=\[2a(\t)+\g_1(t)\](t-\t)^4,\qq\lim_{t\da\t}\g_1(t)=0,\\
\ns\ds\L_2(t)=\g_2(t)(t-\t)^4,\qq\lim_{t\da\t}\g_2(t)=0.\ea$$
Hence,
$$\ba{ll}
\ns\ds\Th(t,\t)\Psi(t,\t)^{-1}={1\over F(t)}\begin{pmatrix}\L_1(t)&\L_2(t)\\ 0&F(t)\end{pmatrix}\\
\ns\ds={1\over\big[2\bar a(\t)+R(t)\big](t-\t)^4}\begin{pmatrix}\big[2a(\t)+\g_1(t)\big](t-\t)^4&\g_2(t)(t-\t)^4\\
0&\big[2\bar a(\t)+R(t)\big](t-\t)^4\end{pmatrix}\\
\ns\ds={1\over2\bar a(\t)+R(t)}\begin{pmatrix}2a(\t)+\g_1(t)&\g_2(t)\\ 0&2\bar a(\t)+R(t)\end{pmatrix}.\ea$$
As a result, we obtain
$$\ba{ll}
\ns\ds\Big|\Th(t,\t)\Psi(t,\t)^{-1}\Big|\les K,\qq\forall t>\t,\hb{ with $t-\t$ small}.\ea$$
The above shows that \eqref{M} holds.

\ex

As we mentioned earlier, condition \eqref{DD>d} implies that $m\ges n$, and if $m=n$ and \eqref{DD>d} holds, then $\h B(\cd)=0$. Hence, in order $\Psi$ to be invertible, one must have $m>n$. The following result is concerned with a case that $D(\cd)$ is surjective and $m=n$.

\bp{m=n} \sl Let {\rm (H1)}, \eqref{DD>d} and
\eqref{Sec2_ABD_Epsilon} hold and $m=n$. Then, for any $p>1$, the
system $[A(\cd),C(\cd);B(\cd),$ $D(\cd)]$ is not $L^p$-exactly
controllable on any $[0,T]$ by $\cU^p[0,T]$ with $T>0$.

\ep

\it Proof. \rm In the current case, $D(\cd)^{-1}$ is bounded. If for
any $x\in\dbR^n$ and $\xi\in L^p_{\cF_T}(\O;\dbR^n)$, one can find a
$u(\cd)\in\cU^p[0,T]$ such that $X(0)=x$ and $X(T)=\xi$, then we let
$$Z(t)=C(t)X(t)+D(t)u(t),\qq t\in[0,T],$$
which will lead to
$$u(t)=D(t)^{-1}\[Z(t)-C(t)X(t)\],\qq t\in[0,T].$$
Hence, $(X(\cd),Z(\cd))$ is an adapted solution to the following BSDE:
$$\left\{\2n\ba{ll}
\ns\ds dX(t)=\[\(A(t)-B(t)D(t)^{-1}C(t)\)X(t)+B(t)D(t)^{-1}Z(t)\]dt+Z(t)dW(t),\q t\in[0,T],\\
\ns\ds X(T)=\xi.\ea\right.$$
Then $X(0)$ cannot be arbitrarily specified. Hence, $L^p$-exact controllability is not possible for system $[A(\cd),C(\cd);B(\cd),D(\cd)]$. \endpf

\ms

In the above two subsections, we have discussed the two extreme cases: either $D(\cd)=0$, or $D(\cd)$ is full rank (for the case $d=1$). The case in between remains open. Some partial results have been obtained, but they are not at a mature level to be reported. We hope to present them in a forthcoming paper.

\section{Duality and Observability Inequality}\label{S_Observability}

As we know that for deterministic linear ODE systems, the
controllability of the original systems is equivalent to the
observability of the dual equations. We would like to see how such a
result will look like for our FSDE system
$[A(\cd),C(\cd);B(\cd),D(\cd)]$. To this end, let us first look at
an abstract result, whose proof should be standard. But for reader's
convenience, we present a proof.

\bp{onto}\sl Let $\dbX$ and $\dbY$ be Banach spaces, $\dbK:\dbX\to\dbY$ be a bounded linear operator, and $\dbK^*:\dbY^*\to\dbX^*$ be the adjoint operator of $\dbK$. Then $\dbK$ is surjective if and only if there exists a $\d>0$ such that
\bel{>d}|\dbK^*y^*|_{\dbX^*}\ges\d|y^*|_{\dbY^*},\qq\forall y^*\in\dbY^*.\ee
Further, if $\dbX$ and $\dbY$ are reflexive and the map $x^*\mapsto|x^*|_{\dbX^*}^2$ from $\dbX^*$ to $\dbR$ is Fr\'echet differentiable, then \eqref{>d} is also equivalent to the following: For any $y\in\dbY$, the functional
\bel{J(y)}J(y^*;y)={1\over2}|\dbK^*y^*|_{\dbX^*}^2+\lan
y,y^*\ran,\qq y^*\in\dbY^*,\ee
admits a minimum over $\dbY^*$. In addition, if the norm of $\dbX^*$ is strictly convex, then for any $y\in\dbY$, the optimal solution of \eqref{J(y)} is necessarily unique.

\ep

\it Proof. \rm Suppose $\cR(\dbK)=\dbY$, i.e., $\mathbb K$ is a
surjection. Then $\cR(\dbK)$ is closed. By Banach Closed Range
Theorem (\cite{Yosida 1980}), $\cR(\dbK^*)$ is closed. Moreover,
$$\cN(\dbK^*)^\perp\equiv\Big\{y\in\dbY\bigm|\lan y^*,y\ran=0,\q\forall y^*\in\cN(\dbK^*)\Big\}
=\cR(\dbK)=\dbY.$$
This implies that $\cN(\dbK^*)=\{0\}$. Thus, $\dbK^*$ is injective
with $\cR(\dbK^*)$ closed. Therefore, $\dbK^*:\dbY^*\to\cR(\dbY^*)$ is
one-to-one and onto. Hence, $(\dbK^*)^{-1}:\cR(\dbK^*)\to\dbY^*$ is
bounded. Consequently, for any $x^*=\dbK^*y^*\in\cR(\dbK^*)$,
$$|y^*|_{\dbY^*}=|(\dbK^*)^{-1}\dbK^*y^*|_{\dbY^*}\les\|(\dbK^*)^{-1}\|
|\dbK^*y^*|_{\dbX^*}\equiv{1\over\d}|\dbK^*y^*|_{\dbX^*},$$
which leads to \eqref{>d}.

\ms

Conversely, suppose (\ref{>d}) holds. Then $\cR(\dbK^*)$ is closed
and $\dbK^*$ is injective. Thus, by Banach Closed Range Theorem,
$\cR(\dbK)$ is also closed and
$$\cR(\dbK)=\cN(\dbK^*)^\perp=\{0\}^\perp=\dbY,$$
proving that $\dbK$ is surjective.

\ms
\ms

Now, for any $y\in\dbY$, consider functional $y^*\mapsto J(y^*;y)$. Clearly, under condition \eqref{>d}, we have that $y^*\mapsto J(y^*;y)$ is coercive and weakly lower semi-continuous. Hence, if $\{y^*_k\}_{k\ges1}$ is a minimizing sequence, then it is bounded. By the reflexivity of $\dbY^*$, we may assume that $y^*_k$ converges weakly to some $\bar y^*\in\dbY^*$. Then by the weakly lower semi-continuity of the functional $J(\cd\,;y)$, $\bar y^*$ must be a minimum.

\ms

Conversely, for any $y\in\dbY$, let $\bar y^*\in\dbY^*$ be a minimum of $y^*\mapsto J(y^*;y)$. Denote the Fr\'echet derivative of $x^*\mapsto{1\over2}|x^*|_{\dbX^*}^2$ by $\G(x^*)$, i.e.,
$$\lim_{\d\to0}{|x^*+\d\xi^*|_{\dbX^*}^2-|x^*|^2_{\dbX^*}\over2\d} \equiv \frac{d}{d\delta}\frac 1 2 | x^* +\delta\xi^* |^2_{\mathbb X^*}\Big|_{\delta=0}  =\lan\G(x^*),\xi^*\ran,\qq
\forall\xi^*\in\dbX^*,$$
where we denote by $\frac{d}{d\delta}\big|_{\delta=0} f(\delta)$ the
derivative of function $f$ at $\delta=0$. Then by the optimality of
$\bar y^*$, we have
\bel{diff}
\ba{ll}
\ns\ds0\les\lim_{\d\to0}{J(\bar y^*+\d y^*;y)-J(\bar y^*;y)\over\d}=\lim_{\d\to0}{|\dbK^*(\bar y^*+\d y^*)|_{\dbX^*}^2-|\dbK^*\bar y^*|_{\dbX^*}^2\over2\d}+\lan y,y^*\ran\\
\ns\ds\q=\lan\G(\dbK^*\bar y^*),\dbK^*y^*\ran+\lan y,y^*\ran=\lan\dbK
\G(\dbK^*\bar y^*)+y,y^*\ran,\qq\forall y^*\in\dbY^*.\ea
\ee
Hence,
$$\dbK\G(\dbK^*\bar y^*)+y=0.$$
Since $y\in\dbY$ is arbitrary, we obtain that $\cR(\dbK)=\dbY$. Then by what we have proved, \eqref{>d} holds. Finally, if we further assume that the norm of $\dbX^*$ is strictly convex, then we must have the uniqueness of the optimal solution to $y^*\mapsto J(y^*;y)$. \endpf

\ms

For any given $1<p<\rho\les\infty$ and $2<\si\les\infty$, we denote
$q\equiv{p\over p-1}$, and
\bel{pqmn}\ba{ll}
\ns\ds\bar p\equiv\left\{\2n\ba{ll}
\ns\ds{p\rho\over\rho-p},\qq\rho<\infty,\\
\ns\ds p,\qq\qq \rho=\infty,\ea\right. \qq\qq
\bar q\equiv\left\{\2n\ba{ll}
\ns\ds{p\rho\over p\rho-\rho+p}, \qq\rho<\infty,\\
\ns\ds{p\over p-1},\qq\qq\rho=\infty,\ea\right.\\
\ns\ds\m\equiv\left\{\2n\ba{ll}
\ns\ds{2\si\over \si-2},\qq \si<\infty,\\
\ns\ds2,\qq\qq \si=\infty,\ea\right.\qq\qq\nu\equiv\left\{\2n\ba{ll}
\ns\ds{2\si\over \si+2},\qq \si<\infty,\\
\ns\ds2,\qq\qq \si=\infty.\ea\right.\ea\ee
Clearly, with the above notations, we have
\bel{=1}\left\{\2n\ba{ll}
\ns\ds{1\over p}+{1\over q}=1,\qq{1\over\bar p}+{1\over\bar q}=1,\qq{1\over\m}+{1\over\n}=1,\\ [4mm]
\ns\ds1<\n\les2\les\m<\infty,\qq1<p\les\bar p<\infty,\qq1<\bar q\les
{p\over p-1},\\ [4mm]
\ns\ds\dbU^{p,\rho,\sigma}[0,T]\equiv L^{\bar
p}_\dbF(\O;L^\m(0,T;\dbR^m)),\qq\dbU^{p,\rho,\si}[0,T]^*\equiv
L^{\bar q}_\dbF(\O;L^\n(0,T;\dbR^m)),\\ [3mm]
\ns\ds\dbU_r^{p,\rho,\si}[0,T]\equiv L^\m_\dbF(0,T;L^{\bar
p}(\O;\dbR^m)),\qq\dbU_r^{p,\rho,\si}[0,T]^*\equiv
L^\n_\dbF(0,T;L^{\bar q}(\O;\dbR^m)).\ea\right.\ee
In the rest of this paper, we will keep the above notations. We now
present the first main result of this section.

\bt{equivalence} \sl Let {\rm(H1)--(H2)} (respectively, {\rm
(H1)} and {\rm(H2)$'$}) hold. Then system $[A(\cd),C(\cd);B(\cd),D(\cd)]$ is
$L^p$-exactly controllable on $[0,T]$ by $\dbU^{p,\rho,\si}[0,T]$
(respectively, $\dbU^{p,\rho,\si}_r[0,T]$) if and only if there
exists a $\d>0$ such that the following, called {\it an
observability inequality}, holds:
\bel{observability inequality}\ba{ll}
\ns\ds\Big\|B(\cd)^\top Y(\cd)+\sum_{k=1}^dD_k(\cd)^\top
Z_k(\cd)\Big\|_{\dbU^{p,\rho,\si}[0,T]^*}\ges\d\|\eta\|_{L^q_{\cF_T}(\O;\dbR^n)},
\qq\forall\eta\in L^q_{\cF_T}(\O;\dbR^n),\ea\ee
(respectively,
\bel{weak observability inequality}\ba{ll}
\ns\ds\Big\|B(\cd)^\top Y(\cd)+\sum_{k=1}^dD_k(\cd)^\top
Z_k(\cd)\Big\|_{\dbU_r^{p,\rho,\si}[0,T]^*}\ges\d\|\eta\|_{L^q_{\cF_T}(\O;\dbR^n)},
\qq\forall\eta\in L^q_{\cF_T}(\O;\dbR^n),\ea\big)\ee
where $(Y(\cd),Z(\cd))$ (with $Z(\cd)\equiv(Z_1(\cd),\cds,Z_d(\cd))$)
is the unique adapted solution to the following BSDE:
\bel{Sec3_Dual_Sys}\left\{\2n\ba{ll}
\ns\ds dY(t)=-\[A(t)^\top Y(t)+\sum_{k=1}^dC_k(t)^\top Z_k(t)\]dt+\sum_{k=1}^dZ_k(t)dW_k(t),\qq t\in [0,T],\\
\ns\ds Y(T)=\eta.\ea\right.\ee

\et

\it Proof. \rm We only prove the equivalence between system's
$L^p$-exactly controllability on $[0,T]$ by $\dbU^{p,\rho,\si}[0,T]$
and the validity of the observability inequality
\eqref{observability inequality}. The other part can be proved with
the similar procedure.

For any $(x,u(\cd))\in\dbR^n\times\dbU^{p,\rho,\si}[0,T]$, with
$p>1$, let $X(\cd)\equiv X(\cd\,;x,u(\cd))$ be the unique solution
to \eqref{FSDE1}. Then
$$X(\cd\,;x,u(\cd))=X(\cd\,;x,0)+X(\cd\,;0,u(\cd)).$$
Define
$$\dbK u(\cd)=X(T;0,u(\cd)),\qq\forall u(\cd)\in\dbU^{p,\rho,\si}[0,T].$$
Then
$$\|\dbK u(\cd)\|_{L^p_{\cF_T}(\O;\dbR^n)}\les K\|u(\cd)\|_{\dbU^{p,\rho,\si}[0,T]},\qq u(\cd)\in\dbU^{p,\rho,\si}[0,T].$$
Thus, $\dbK:\dbU^{p,\rho,\si}[0,T]\to L^p_{\cF_T}(\O;\dbR^n)$ is a
bounded linear operator. Now, system $[A(\cd),C(\cd);B(\cd),D(\cd)]$
is $L^p$-exactly controllable on $[0,T]$ by $\dbU^{p,\rho,\si}[0,T]$
if and only if for any $(x,\xi)\in\dbR^n\times
L^p_{\cF_T}(\O;\dbR^n)$,
$$\xi-X(T;x,0)=\dbK u(\cd),$$
for some $u(\cd)\in\dbU^{p,\rho,\si}[0,T]$, which is equivalent to
$$\cR(\dbK)=L^p_{\cF_T}(\O;\dbR^n).$$
This means that $\dbK:\dbU^{p,\rho,\si}[0,T]\to
L^p_{\cF_T}(\O;\dbR^n)$ is surjective. Hence, by Proposition
\ref{onto}, this is equivalent to the following: For some $\d>0$,
\bel{|K*eta|}\|\dbK^*\eta\|_{\dbU^{p,\rho,\si}[0,T]^*}\ges\d\|\eta\|_{L^q_{\cF_T}(\O:\dbR^n)},
\qq\forall\eta\in L^q_{\cF_T}(\O;\dbR^n).\ee
We now find
$\dbK^*:L^q_{\cF_T}(\O;\dbR^n)\to\dbU^{p,\rho,\si}[0,T]^*$. To this
end, we consider BSDE \eqref{Sec3_Dual_Sys}, with $\eta\in
L^q_{\cF_T}(\O;\dbR^n)$. By a standard result of BSDEs (\cite{El
Karoui-Peng-Quenez 1997}), under (H1), \eqref{Sec3_Dual_Sys} admits
a unique adapted solution
$$(Y(\cd),Z(\cd))\equiv\big(Y(\cd\,;\eta),Z(\cd\,;\eta)\big)\in L^q_\dbF(\O;C([0,T];\dbR^n))\times L^q_\dbF(\O;L^2(0,T;\dbR^{n\times d})),$$
and the following estimate holds:
\bel{}\dbE\[\sup_{t\in[0,T]}|Y(t)|^{q}+\(\int_0^T|Z(t)|^2dt\)^{q\over2}\]\les
K\dbE|\eta|^{q},\ee
where we denote $Z(\cd)\equiv(Z_1(\cd),Z_2(\cd),\dots,Z_d(\cd))$. By applying
It\^o's formula to $\lan X(\cd\,;0,u(\cd)),Y(\cd)\ran$ on the interval $[0,T]$, we have the following duality relation:
\bel{Sec3_Duality}\ba{ll}
\ns\ds\lan\dbK u(\cd),\eta\ran=\dbE\lan X(T),\eta\ran=\dbE\int_0^T\lan u(t),B(t)^\top Y(t)+\sum_{k=1}^dD_k(t)^\top Z_k(t)\ran dt=\lan u(\cd),\dbK^*\eta\ran.\ea\ee
Hence,
\bel{K*}(\dbK^*\eta)(t)\equiv B(t)^\top Y(t)+\sum_{k=1}^dD_k(t)^\top
Z_k(t),\quad t\in [0,T].\ee
Combining \eqref{|K*eta|} with \eqref{K*}, we obtain \eqref{observability inequality}. \endpf

\ms

Theorem \ref{equivalence} provides an approach to study the controllability of stochastic linear systems by establishing an inequality for BSDEs. The following example illustrates this approach.

\bex{Example 4.3} \rm Let the dimensions of both state process and Brownian motion be 1, the dimension of control process be 2. Let
$$A(\cd), B_1(\cd), D_2(\cd)\in L^\infty_\dbF(0,T;\dbR),\q| B_1(t)|,\ | D_2(t)|\ges\d,\qq t\in[0,T],$$
for some $\d>0$. Consider the following system:
\bel{Sec3_Example_Orig_Sys}\left\{\2n\ba{ll}
\ns\ds dX(t)=\Bigg[A(t)X(t)+\begin{pmatrix} B_1(t) & 0\end{pmatrix}\begin{pmatrix}u_1(t)\\ u_2(t)\end{pmatrix}\Bigg]dt+\begin{pmatrix} 0 & D_2(t)\end{pmatrix}\begin{pmatrix}u_1(t)\\ u_2(t)\end{pmatrix}dW(t),\q t\ges 0,\\
\ns\ds X(0)=x,\ea\right.\ee
and the adjoint system is given by
\bel{Sec3_Example_Dual_Sys}\left\{\2n\ba{ll}
\ns\ds dY(t)=-A(t)Y(t)dt+Z(t)dW(t),\qq t\in[0,T],\\
\ns\ds y(T)=\eta,\ea\right.\ee
where $\eta\in L^2_{\mathbb F}(0,T;\mathbb R)$. A direct calculation leads to
\bel{Sec3_Example_Eq_0}\ba{ll}
\ns\ds\dbE\int_0^T\Bigg|\begin{pmatrix} B_1(t)\\ 0\end{pmatrix}Y(t)+\begin{pmatrix}0\\  D_2(t)\end{pmatrix}Z(t)\Bigg|^2dt=\dbE\int_0^T\[| B_1(t)|^2|Y(t)|^2
+| D_2(t)|^2|Z(t)|^2\]dt\\ [3mm]
\ns\ds\qq\qq\qq\qq\qq\qq\qq\qq\qq\ges\d^2\dbE\int_0^T\[|Y(t)|^2+|Z(t)|^2\]dt.\ea\ee
Let $\b\in L^\infty_\dbF(0,T;\dbR)$ to be chosen later. Applying It\^o's formula to
$|Y(\cd)|^2e^{\int_0^\cd\b(s)ds}$ on the interval $[0,T]$, we deduce that
\bel{Sec3_Example_Eq_1}\dbE\[|\eta|^2e^{\int_0^T\b(s)ds}\]-|Y(0)|^2=\dbE\int_0^T\[\big(\b(t)-2A(t)\big)
|Y(t)|^2+|Z(t)|^2\]e^{\int_0^t\b(s)ds}dt.\ee
In \eqref{Sec3_Example_Eq_1}, selecting $\b(\cd)=2A(\cd)+1$ leads to
\bel{Sec3_Example_Eq_2}\dbE\int_0^T\[|Y(t)|^2+|Z(t)|^2\]e^{\int_0^t
(2A(s)+1)ds}dt=\dbE\[|\eta|^2e^{\int_0^T(2A(s)+1)ds}\]-|Y(0)|^2.\ee
On the other hand, selecting $\b(\cd)=2A(\cd)$, we get
$$\dbE\[|\eta|^2e^{\int_0^T2A(s)ds}\]-|Y(0)|^2=\dbE\int_0^T|Z(t)|^2e^{\int_0^t2A(s)ds}dt\ges0.$$
Then,
\bel{Sec3_Example_Eq_3}-|Y(0)|^2\ges-\dbE\[|\eta|^2e^{\int_0^T2A(s)ds}\].\ee
Combining \eqref{Sec3_Example_Eq_2} with \eqref{Sec3_Example_Eq_3}, one has
$$\ba{ll}
\ns\ds\dbE\int_0^T\[|Y(t)|^2+|Z(t)|^2\]e^{\int_0^t(2A(s)+1)ds}dt\ges\dbE\[|\eta|^2e^{\int_0^T
(2A(s)+1)ds}\]-\dbE\[|\eta|^2e^{\int_0^T2A(s)ds}\]\\
\ns\ds\qq\qq\qq\qq\qq\qq\qq\qq\q~=(e^T-1)\dbE\[|\eta|^2e^{\int_0^T2A(s)ds}\].\ea$$
Due to the boundedness of $A$ (without loss of generality, we assume
$|A(t)|\les K$ a.s. a.e.), we obtain
$$e^{(2K+1)T}\dbE\int_0^T\[|Y(t)|^2+|Z(t)|^2\]dt\ges(e^T-1)e^{-2KT}\dbE|\eta|^2,$$
i.e.
$$\dbE\int_0^T\[|Y(t)|^2+|Z(t)|^2\]dt\ges(e^T-1)e^{-(4K+1)T}\dbE|\eta|^2.$$
Combining with \eqref{Sec3_Example_Eq_0}, we have proved that the
observability inequality holds true for BSDE
\eqref{Sec3_Example_Dual_Sys}. By Theorem \ref{equivalence}, system
\eqref{Sec3_Example_Orig_Sys} is $L^2$-exactly controllable on
$[0,T]$ by $\dbU^{2,\infty,\infty}[0,T]\equiv L^2_{\mathbb
F}(0,T;\mathbb R^2)$.

\ex

Now, we introduce the following definition which makes the name ``observability inequality'' aforementioned meaningful.

\bde{observability} \rm Let (H1) hold and $(Y(\cd),Z(\cd))$ be the adapted solution to BSDE \eqref{Sec3_Dual_Sys} with $\eta\in L_{\cF_T}^q(\O;\dbR^n)$.

\ms

(i) For the pair $(B(\cd),D(\cd))$ with $B(\cd),D_k(\cd)\in L^1_\dbF(0,T;\dbR^{n\times m})$ ($k=1,2,\cds,d$) and $D(\cd)=(D_1(\cd),\cds,$ $D_d(\cd))$, the map
\bel{observer}\eta\mapsto\dbK^*\eta\equiv B(\cd)^\top Y(\cd)+\sum_{k=1}^dD_k(\cd)^\top Z_k(\cd)\ee
is called an $\dbY[0,T]$-{\it observer} of BSDE \eqref{Sec3_Dual_Sys} if
$$\ds \dbK^*\eta\in\dbY[0,T],\qq\forall\eta\in L^q_{\cF_T}(\O;\dbR^m),$$
where $\dbY[0,T]$ is a subspace of $L^1_\dbF(0,T;\dbR^m)$. BSDE \eqref{Sec3_Dual_Sys}, together with the observer \eqref{observer} is denoted by $[A(\cd)^\top,C(\cd)^\top;B(\cd)^\top,D(\cd)^\top]$.

\ms

(ii) System $[A(\cd)^\top,C(\cd)^\top;B(\cd)^\top,D(\cd)^\top]$ is said to be {\it $L^q$-exactly observable} by $\dbY[0,T]$ observations if from the {\it observation} $\dbK^*\eta\in\dbY[0,T]$, the terminal value $\eta\in L^q_{\cF_T}(\O;\dbR^n)$ of $Y(\cd)$ at $T$ can be uniquely determined, i.e., the map $\dbK^*:L^q_{\cF_T}(\O;\dbR^n)\to\dbY[0,T]$ admits a bounded inverse.

\ede

With the above definition, we clearly have the following result:

\bt{observe} \sl
Let {\rm (H1)--(H2)} (respectively, {\rm (H1)--(H2)$'$}) hold. Then $[A(\cd),C(\cd);B(\cd),D(\cd)]$ is $L^p$-exactly controllable on $[0,T]$ by $\dbU^{p,\rho,\si}[0,T]$ (respectively, $\dbU_r^{p,\rho,\si}[0,T]$) if and only if $[A(\cd)^\top,C(\cd)^\top;B(\cd)^\top,D(\cd)^\top]$ is $L^q$-exactly
observable by $\dbU^{p,\rho,\si}[0,T]^*$ (respectively,
$\dbU_r^{p,\rho,\si}[0,T]^*$) observations.

\et

Next, for any $x\in\dbR^n$, let $X(\cd\,;x,0)$ be the solution to the state equation \eqref{FSDE1} corresponding to the initial state $x$ and $u(\cd)=0$. Denote
$$\dbK_0x=X(T;x,0),\qq\forall x\in\dbR^n.$$
Then applying It\^o's formula to $\lan X(\cd\,;x,0),Y(\cd)\ran$ with $(Y(\cd),Z(\cd))\equiv(Y(\cd\,;\eta),Z(\cd\,;\eta))$ being the adapted solution to \eqref{Sec3_Dual_Sys}, we have
$$\ba{ll}
\ns\ds\dbE\lan X(T),\eta\ran-\lan x,Y(0)\ran=\dbE\int_0^T\(\lan A(t)X(t),Y(t)\ran+\lan X(t),-A(t)^\top Y(t)-\sum_{k=1}^dC_k(t)^\top Z_k(t)\ran\\
\ns\ds\qq\qq\qq\qq\qq\qq\qq+\sum_{k=1}^d\lan C_k(t)X(t),Z_k(t)\ran\)dt=0.\ea$$
Hence,
$$\lan x,Y(0)\ran=\dbE\lan\dbK_0x,\eta\ran=\lan x,\dbK_0^*\eta\ran,\qq\forall x\in\dbR^n.$$
This leads to
\bel{K0*}\dbK_0^*\eta=Y(0;\eta),\qq\forall\eta\in L^q_{\cF_T}(\O;\dbR^n).\ee

\ms

Now, for any $(x,\xi)\in\dbR^n\times L^p_{\cF_T}(\O;\dbR^n)$, we introduce a functional $J(\cd\,;x,\xi):L^q_{\cF_T}(\O;\dbR^n)\to\dbR$ defined by
\bel{Sec3_Dual_Cost}J(\eta;x,\xi)={1\over2}\|\dbK^*\eta\|_{\dbU^{p,\rho,\si}[0,T]^*}^2+\lan
x,\dbK_0^*\eta\ran-\dbE\lan\xi,\eta \ran,\qq\forall\eta\in
L^q_{\cF_T}(\O;\dbR^n),\ee
where $\dbK^*$ and $\dbK_0^*$ are given by \eqref{K*} and \eqref{K0*}, respectively. Equivalently,
\bel{Sec3_Dual_Cost*}\ba{ll}
\ns\ds J(\eta;x,\xi)={1\over2}\|B(\cd)^\top
Y(\cd)+\sum_{k=1}^dD_k(\cd)^\top
Z_k(\cd)\|_{\dbU^{p,\rho,\si}[0,T]^*}^2+\lan x,Y(0)\ran
-\dbE\lan\xi,\eta\ran,\q\forall\eta\in L^q_{\cF_T}(\O;\dbR^n),\ea\ee
with $(Y(\cd),Z(\cd))$ being the adapted solution to BSDE \eqref{Sec3_Dual_Sys}. One can pose the following optimization problem.

\ms

\bf Problem (O). \rm Minimize \eqref{Sec3_Dual_Cost*} subject to BSDE
\eqref{Sec3_Dual_Sys} over $L^q_{\cF_T}(\O;\dbR^n)$.

\ms

Note that the spaces $\dbU^{p,\rho,\si}[0,T]$ and
$L^q_{\cF_T}(\O;\dbR^n)$ are reflexive since their norms are
uniformly convex. In order to apply Proposition \ref{onto}, we need
to show that
$$\f(\cd)\mapsto{1\over2}\|\f(\cd)\|^2_{\dbU^{p,\rho,\si}[0,T]^*}\equiv{1\over2}
\|\f(\cd)\|^2_{L^{\bar q}_\dbF(\O;L^\n(0,T;\dbR^m))}$$
is Fr\'echet differentiable. For simplicity of notation, we shall
use the following notations for a while:
$$
\| \cdot \|_{L^{\bar q}_{\mathbb F} (\Omega;L^\nu(0,T; \mathbb
R^m))} \equiv \| \cdot \|_{L^{\bar q,\nu}_{\mathbb F}} \mbox{ and }
\|\cdot\|_{L^\nu(0,T; \mathbb R^m)}\equiv \| \cdot \|_{L^\nu}.
$$
Now let $\varphi(\cdot), \psi(\cdot)\in L^{\bar q}_{\mathbb
F}(\Omega;L^\nu(0,T;\mathbb R^m))$. If $\|\varphi(\cdot)\|_{L^{\bar
q,\nu}_{\mathbb F}}=0$, then
\begin{equation}\label{Sec3_Eq_Temp1.1}
\frac 1 2 \frac{d}{d\delta}\Big\{ \|
\varphi(\cdot)+\delta\psi(\cdot) \|_{L^{\bar q,\nu}_{\mathbb F}}^2
\Big\}\Big|_{\delta =0}  =0.
\end{equation}
If $\|\varphi(\cdot)\|_{L^{\bar q,\nu}_{\mathbb F}} \neq 0$,
\begin{equation}\label{Sec3_Eq_Temp1.2}
\begin{aligned}
\frac 1 2 \frac{d}{d\delta}\Big\{ \|
\varphi(\cdot)+\delta\psi(\cdot) \|_{L^{\bar q,\nu}_{\mathbb F}}^2
\Big\}\Big|_{\delta =0} =\ & \frac 1 2 \frac{d}{d\delta}\bigg\{
\bigg[ \mathbb E\bigg( \int_0^T |\varphi(t) +\delta\psi(t)|^{\nu} dt
\bigg)^{\frac{\bar q}{\nu}} \bigg]^{\frac {2}{\bar q}}
\bigg\}\Big|_{\delta =0}\\
=\ & \frac {1}{\bar q} \| \varphi(\cdot) \|_{L^{\bar q,\nu}_{\mathbb
F}}^{2-\bar q} \frac{d}{d\delta}\Big\{\mathbb
E[f(\omega,\delta)] \Big\}\Big|_{\delta =0},\\
\end{aligned}
\end{equation}
provided the derivative on the right hand side exists, where
$$f(\omega,\delta) = \bigg( \int_0^T |\varphi(t,\omega) +\delta
\psi(t,\omega)|^{\nu} dt \bigg)^\frac{\bar q}{\nu}.$$
For simplicity of notation, we denote $\infty\times 0 =0$. Then
\eqref{Sec3_Eq_Temp1.1} and \eqref{Sec3_Eq_Temp1.2} can be combined
into
\begin{equation}\label{Sec3_Eq_Temp2}
\frac 1 2 \frac{d}{d\delta}\Big\{ \|\varphi(\cdot)
+\delta \psi(\cdot)\|_{L^{\bar q,\nu}_\dbF}^2 \Big\}\Big|_{\delta =0} = \frac
{1}{\bar q} \| \varphi(\cdot) \|_{L^{\bar q,\nu}_{\mathbb
F}}^{2-\bar q} {\bf1}_{\{ \|\varphi(\cdot)\|_{L^{\bar q,\nu}_\dbF}
\neq 0 \}} \frac{d}{d\delta}\Big\{\mathbb
E[f(\omega,\delta)] \Big\}\Big|_{\delta =0}.
\end{equation}
To exchange the order of derivation and expectation in
\eqref{Sec3_Eq_Temp2}, we calculate $\frac{\partial
f}{\partial\delta}$. For any $\delta\in (-1,1)$ and
$\omega\in\Omega$, if $\| \varphi(\cdot,\omega)
+\delta\psi(\cdot,\omega) \|_{L^\nu} =0$, i.e., $\varphi(t,\omega)
+\delta\psi(t,\omega) =0$ a.e. $t\in [0,T]$, then
\begin{equation}\label{sec3_Eq_Temp2.1}
\frac{\partial f}{\partial\delta} (\omega,\delta) =
\lim_{\Delta\delta\rightarrow 0} \frac{1}{\Delta\delta} \Big(
f(\omega,\delta+\Delta\delta) -f(\omega,\delta) \Big) =
\lim_{\Delta\delta\rightarrow 0} \frac{1}{\Delta\delta}
\|\Delta\delta\psi(\cdot,\omega)\|_{L^\nu}^{\bar q} =0.
\end{equation}
On the other hand, when $\| \varphi(\cdot,\omega)
+\delta\psi(\cdot,\omega) \|_{L^\nu} \neq 0$,
\begin{equation}\label{sec3_Eq_Temp2.2}
\frac{\partial f}{\partial\delta}(\omega,\delta) = \frac{\bar
q}{\nu} \|\varphi(\cdot,\omega)
+\delta\psi(\cdot,\omega)\|_{L^{\nu}}^{\bar q-\nu}
\frac{\partial}{\partial \delta} \bigg(\int_0^T g(t,\omega,\delta)
dt\bigg),
\end{equation}
where
$$
g(t,\omega,\delta) = |\varphi
(t,\omega)+\delta\psi(t,\omega)|^{\nu}.
$$
Combining \eqref{sec3_Eq_Temp2.1} with \eqref{sec3_Eq_Temp2.2}, one
has
\begin{equation}\label{sec3_Eq_Temp2.3}
\frac{\partial f}{\partial\delta}(\omega,\delta) = \frac{\bar
q}{\nu} \|\varphi(\cdot,\omega)
+\delta\psi(\cdot,\omega)\|_{L^{\nu}}^{\bar q-\nu} {\bf1}_{\{
\|\varphi(\cdot,\omega) +\delta\psi(\cdot,\omega)\|_{L^{\nu}} \neq 0
\}} \frac{\partial}{\partial \delta} \bigg(\int_0^T
g(t,\omega,\delta) dt\bigg).
\end{equation}
To exchange the order of derivation and integral in
\eqref{sec3_Eq_Temp2.3}, we calculate
$$
\frac{\partial g}{\partial\delta}(t,\omega,\delta) = \nu
|\varphi(t,\omega) +\delta\psi(t,\omega)|^{\nu-2} \langle
\varphi(t,\omega) +\delta\psi(t,\omega),\ \psi(t,\omega) \rangle,
$$
and then
$$
\begin{aligned}
\bigg| \frac{\partial g}{\partial\delta}(t,\omega,\delta) \bigg|
\les\ & \nu |\varphi(t,\omega) +\delta\psi(t,\omega)|^{\nu-1}
|\psi(t,\omega)|\\
\les\ & C \Big( |\varphi(t,\omega)|^\nu +|\psi(t,\omega)|^\nu \Big)
\in L^1(0,T;\mathbb R),
\end{aligned}
$$
where $C>0$ is a constant only depending on $\nu$. By Theorem 2.27
in \cite[Page 56]{Folland 1999}, the order of derivation and
integral in \eqref{sec3_Eq_Temp2.3} can be exchanged. Then, we have
\begin{equation}\label{sec3_Eq_Temp2.4}
\begin{aligned}
\frac{\partial f}{\partial\delta}(\omega,\delta) =&\ \frac{\bar
q}{\nu} \|\varphi(\cdot,\omega)
+\delta\psi(\cdot,\omega)\|_{L^{\nu}}^{\bar q-\nu}{\bf1}_{\{
\|\varphi(\cdot,\omega)
+\delta\psi(\cdot,\omega)\|_{L^{\nu}} \neq 0 \}} \int_0^T \frac{\partial g}{\partial \delta} (t,\omega,\delta) dt\\
=&\  \bar q \|\varphi(\cdot,\omega)
+\delta\psi(\cdot,\omega)\|_{L^{\nu}}^{\bar q-\nu}{\bf1}_{\{
\|\varphi(\cdot,\omega) +\delta\psi(\cdot,\omega)\|_{L^{\nu}} \neq 0
\}}\\
& \times\int_0^T |\varphi
(t,\omega)+\delta\psi(t,\omega)|^{\nu-2}\langle
\varphi(t,\omega)+\delta\psi(t,\omega),\ \psi(t,\omega) \rangle dt\\
=&\  \bar q \|\varphi(\cdot,\omega)
+\delta\psi(\cdot,\omega)\|_{L^{\nu}}^{\bar q-\nu}\\
& \times\int_0^T |\varphi
(t,\omega)+\delta\psi(t,\omega)|^{\nu-2}\langle
\varphi(t,\omega)+\delta\psi(t,\omega),\ \psi(t,\omega) \rangle dt.
\end{aligned}
\end{equation}
By virtue of H\"{o}lder's inequality, we have
$$\begin{aligned}
\bigg|\frac{\partial f}{\partial \delta}(\omega,\delta)\bigg| \les&\
\bar q \|\varphi(\cdot,\omega)
+\delta\psi(\cdot,\omega)\|_{L^{\nu}}^{\bar q-\nu} \int_0^T  |\varphi (t,\omega)+\delta\psi(t,\omega)|^{\nu-1}|\psi(t,\omega)| dt\\
\les&\ \bar q \| \varphi(\cdot,\omega) +\delta\psi(\cdot,\omega)
\|_{L^\nu}^{\bar q-\nu} \|
\varphi(\cdot,\omega)+\delta\psi(\cdot,\omega)
\|_{L^\nu}^{\nu-1} \|\psi(\cdot,\omega)\|_{L^\nu} \\
\les&\ C\left( \|\varphi(\cdot,\omega)\|_{L^\nu}^{\bar q}
+\|\psi(\cdot,\omega)\|_{L^\nu}^{\bar q} \right) \in L^1_{\mathcal
F_T}(\Omega;\mathbb R),
\end{aligned}$$
where $C>0$ is a constant only depending on $\bar q$. Theorem 2.27
in \cite[Page 56]{Folland 1999} works again to exchange the order of
derivation and expectation in \eqref{Sec3_Eq_Temp2}. Then, combining
with \eqref{sec3_Eq_Temp2.4}, we have
\begin{equation}\label{sec3_Eq_Temp2.5}
\begin{aligned}
& \frac 1 2 \frac{d}{d\delta} \Big\{
\|\varphi(\cdot) +\delta \psi(\cdot)\|_{L^{\bar q,\nu}_\dbF}^2
\Big\}\Big|_{\delta =0} = \frac {1}{\bar q} \|\varphi(\cdot)\|_{L^{\bar
q,\nu}_{\mathbb F}}^{2-\bar q} {\bf1}_{\{
\|\varphi(\cdot)\|_{L^{\bar q,\nu}_\dbF} \neq 0 \}} \mathbb E
\bigg\{ \frac{\partial}{\partial\delta}
f(\omega,\delta) \bigg\}\Big|_{\delta =0}\\
=\ & \|\varphi(\cdot)\|_{L^{\bar q,\nu}_\dbF}^{2-\bar q} {\bf1}_{\{
\|\varphi(\cdot)\|_{L^{\bar q,\nu}_\dbF} \neq 0 \}} \mathbb E
\bigg\{ \|\varphi(\cdot,\omega)\|_{L^{\nu}}^{\bar q-\nu} \int_0^T
|\varphi (t)|^{\nu-2}\langle
\varphi(t),\ \psi(t) \rangle dt \bigg\}\\
=\ & \|\varphi(\cdot)\|_{L^{\bar q,\nu}_\dbF}^{2-\bar q} \mathbb E
\bigg\{ \|\varphi(\cdot,\omega)\|_{L^{\nu}}^{\bar q-\nu} \int_0^T
|\varphi (t)|^{\nu-2}\langle
\varphi(t),\ \psi(t) \rangle dt \bigg\} \\
=\ & \mathbb E\int_0^T \bigg\langle \|\varphi(\cdot)\|^{2-\bar
q}_{L^{\bar q,\nu}_{\mathbb F}} \mathbb E_t \Big[ \|
\varphi(\cdot,\omega) \|^{\bar q-\nu}_{L^\nu} {\bf1}_{\{
\|\varphi(\cdot,\omega)\|_{L^\nu}\neq 0 \}} \Big]
|\varphi(t)|^{\nu-2}\varphi(t),\ \psi(t) \bigg\rangle dt\\
\equiv\ & \mathbb E\int_0^T \langle \Gamma(\varphi(\cdot))(t),\
\psi(t) \rangle dt,
\end{aligned}
\end{equation}
where
\begin{equation}\label{G(f)*}
\Gamma(\varphi(\cdot))(t) = \| \varphi(\cdot) \|^{2-\bar q}_{L^{\bar
q}_{\mathbb F}(\Omega;L^\nu(0,T;\mathbb R^m))} M(t)
|\varphi(t)|^{\nu-2} \varphi(t),\quad t\in [0,T],
\end{equation}
with
\bel{M(t)}M(t)=\dbE_t\[ \|\varphi(\cdot,\omega)\|^{\bar
q-\nu}_{L^\nu(0,T;\mathbb R^m)} {\bf1}_{\{
\|\varphi(\cdot,\omega)\|_{L^\nu(0,T;\mathbb R^m)} \neq 0 \}}
\],\qq t\in[0,T].\ee
We have the following lemma.

\bl{Lemma 4.6} \sl Let $p\les \frac{2\si\rho}{\si\rho-2\rho+2\si}$.
Then for any $\f(\cd)\in L^{\bar
q}_\dbF(\O;L^\n(0,T;\dbR^m))\equiv\dbU^{p,\rho,\si}[0,T]^*$,
\bel{}\G(\f(\cd))(\cd)\in L^{\bar
p}_\dbF(\O;L^\m(0,T;\dbR^m))=\dbU^{p,\rho,\si}[0,T].\ee

\el

We notice that, in $p\les \frac{2\si\rho}{\si\rho-2\rho+2\si}$,
$\rho$ takes values in $(1,\infty]$ and $\si$ takes values in
$(2,\infty]$; When $\rho=\infty$ or/and $\si=\infty$, the right hand
side of the inequality takes the limit. Moreover the inequality
$p\les \frac{2\si\rho}{\si\rho-2\rho+2\si}$ is equivalent to
$\bar q \ges \nu$ or $\bar p \les \mu$.\\

\it Proof of Lemma \ref{Lemma 4.6}. \rm Note that $p\les
\frac{2\si\rho}{\si\rho-2\rho+2\si}$ is equivalent to  $\bar
q\ges\n$. It suffices to show that
$$\h\G(\f(\cd))\equiv M(\cd)|\f(\cd)|^{\n-2}\f(\cd)\in L^{\bar q}_\dbF(\O;L^\n(0,T;\dbR^m)).$$
First of all, if $\bar q=\n$, then $\bar p=\m$. In this case, $M(\cd)=1$ and
$$\h\G(\f(\cd))(t)=|\f(t)|^{\n-2}\f(t),\qq t\in[0,T].$$
Hence, (note that $L^{\bar p}_\dbF(\O;L^\m(0,T;\dbR^m))=L^\m_\dbF(0,T;\dbR^m)$ in the current case)
$$\dbE\int_0^T|\h\G(\f(\cd))(t)|^\m dt=\dbE\int_0^T\(|\f(t)|^{\n-1}\)^\m dt=\dbE\int_0^T|\f(t)|^\n dt<\infty.$$
Now, let $\bar q>\n$. Then
$$\m={\n\over\n-1}>{\bar q\over\bar q-1}=\bar p,$$
which leads to ${\bar q\m\over\bar p\n}>1$. Note that $M(\cd)$ is a (nonnegative valued) martingale. By Jensen's inequality,
$$\dbE M(t)^{\bar q\over\bar q-\n}=\dbE\Big\{\dbE_t\[\(\int_0^T|\f(t)|^\n dt\)^{{\bar q\over\n}-1}\]\Big\}^{\bar q\over\bar q-\n}\les\dbE\(\int_0^T|\f(t)|^\n dt\)^{\bar q\over\n}.$$
Hence, using Doob's inequality,
$$\dbE\[\sup_{t\in[0,T]}M(t)^{\bar q\over\bar q-\n}\]\les\({\bar q\over\n}\)^{\bar q\over\bar q-\n}\dbE\big[ M(T)^{\bar q\over\bar q-\n}\big]=\({\bar q\over\n}\)^{\bar q\over\bar q-\n}\dbE\(\int_0^T|\f(t)|^\n dt\)^{\bar q\over\n}.$$
Consequently,
$$\ba{ll}
\ns\ds\dbE\(\int_0^T|\h\G(\f(\cd))(t)|^\m dt\)^{\bar p\over\m}=\dbE\(\int_0^TM(t)^\m|\f(t)|^{(\n-1)\m}dt\)^{\bar p\over\m}\\
\ns\ds\les\dbE\[\sup_{t\in[0,T]}M(t)^{\bar p}\(\int_0^T|\f(t)|^\n dt\)^{\bar p\over\m}\]\\
\ns\ds\les\[\dbE\(\sup_{t\in[0,T]}M(t)^{\bar p\bar q\m\over\bar q\n-\bar p\m}\)\]^{\bar q\m-\bar p\n\over\bar q\m}\[\dbE\(\int_0^T|\f(t)|^\n dt\)^{\bar q\over\n}\]^{\bar p\n\over\bar q\m}.\ea$$
Since
$${\bar p\m\over\bar q\m-\bar p\n}={{\bar q\over\bar q-1}{\n\over\n-1}\over
{\bar q\n\over\n-1}-{\bar q\n\over\bar q-1}}={1\over\bar q-\n},\qq{\bar p\n\over\bar q\m}={\n-1\over\bar q-1},$$
from the above, we obtain
\bel{doob}
\ba{ll}
\ns\ds\dbE\(\int_0^T|\h\G(\f(\cd))(t)|^\m dt\)^{\bar p\over\m}\les\[\dbE\(\sup_{t\in[0,T]}M(t)^{\bar q\over\bar q-\n}\)\]^{\bar p(\bar q-\n)\over\bar q}\[\dbE\(\int_0^T|\f(t)|^\n dt\)^{\bar q\over\n}\]^{\n-1\over\bar q-1}\\
\ns\ds\les\Big\{\({\bar q\over\n}\)^{\bar q\over\bar q-\n}\dbE\(\int_0^T|\f(t)|^\n dt\)^{\bar q\over\n}\Big\}^{\bar p(\bar q-\n)\over\bar q}\[\dbE\(\int_0^T|\f(t)|^\n dt\)^{\bar q\over\n}\]^{\n-1\over\bar q-1}\\
\ns\ds=\({\bar q\over\n}\)^{\bar q\over\bar q-1}\dbE\(\int_0^T|\f(t)|^\n dt\)^{\bar q\over\n}<\infty.\ea\ee
This proves out conclusion. \endpf

\ms

Now, let us look at the optimal solution $\bar\eta$ of Problem (O). According to the above, we see that the optimal solution $\bar\eta$ of Problem (O) satisfies the following:
$$0=\mathbb E \lan\dbK\G(\dbK^*\bar\eta)+\dbK_0x-\xi,\eta\ran,\qq
\forall\eta\in L^q_{\cF_T}(\O;\dbR^n),$$
where $\G(\dbK^*\bar\eta)$ is given by \eqref{G(f)*} with $\f(\cd)=\dbK^*\bar\eta$. Thus,
\bel{4.25}\dbK\G(\dbK^*\bar\eta)+\dbK_0x-\xi=0.\ee
Now, when $p\les \frac{2\si\rho}{\si\rho-2\rho+2\si}$, we define
\bel{bar u}\bar u(\cd)=\G(\dbK^*\bar\eta)\in L^{\bar
p}_\dbF(\O;L^\m(0,T;\dbR^m))\equiv\dbU^{p,\rho,\si}[0,T].\ee
Then \eqref{4.25} reads
$$\xi=\dbK\bar u(\cd)+\dbK_0x=X(T;x,\bar u(\cd)),$$
which means that $\bar u(\cd)\in\dbU^{p,\rho,\si}[0,T]$ is a control
steering $x\in\dbR^n$ to $\xi\in L^p_{\cF_T}(\O;\dbR^n)$. Therefore,
we obtain the following result, making use of Proposition \ref{onto}
and Theorem \ref{equivalence}.

\bt{Theorem 4.7} \sl Let {\rm(H1)--(H2)} hold and $p\les
\frac{2\si\rho}{\si\rho-2\rho+2\si}$. Then the observability
inequality \eqref{observability inequality} holds if and only if for
any $(x,\xi)\in\dbR^n\times L^p_{\cF_T}(\O;\dbR^n)$, Problem
{\rm(O)} admits a unique optimal solution $\bar\eta\in
L^q_{\cF_T}(\O;\dbR^n)$. In this case, the control $\bar
u(\cd)\in\dbU^{p,\rho,\si}[0,T]$ defined by \eqref{bar u} steers $x$
to $\xi$. Moreover, with $\bar u(\cd)$ defined by \eqref{bar u} for
$$\dbK^*\bar\eta=B(\cd)^\top\bar Y(\cd)+\sum_{k=1}^dD_k(\cd)^\top\bar Z_k(\cd),$$
the following coupled FBSDE
\bel{Sec3_Haml_Sys}\left\{\2n\ba{ll}
\ns\ds d\bar X(t)=\[A(t)\bar X(t)+B(t)\bar u(t)\]dt+\sum_{k=1}^d
\[C_k(t)\bar X(t)+D_k(t)\bar u(t)\]dW_k(t),\\
\ns\ds d\bar Y(t)=-\[A(t)^\top\bar Y(t)+\sum_{k=1}^d
C_k(t)^\top\bar Z_k(t)\]dt+\sum_{k=1}^d\bar Z_k(t)dW_k(t),\\
\ns\ds\bar X(0)=x,\qq\bar X(T)=\xi\ea\right.\ee
admits a unique adapted solution
$$(\bar X(\cd),\bar Y(\cd),\bar Z(\cd))\in L^p_\dbF(\O;C([0,T];\dbR^n))\times L^q_\dbF(\O;C([0,T];\dbR^n))\times L^q_\dbF(\O;L^2(0,T;\dbR^{n\times d})).$$
\et

\it Proof. \rm By the above analysis, the only remaining thing is to
prove the uniqueness of FBSDE \eqref{Sec3_Haml_Sys}. Now, let
$(\tilde X(\cdot),\tilde Y(\cdot),\tilde Z(\cdot)) \in L^p_{\mathbb
F}(\Omega;C([0,T];\mathbb R^n)) \times L^q_{\mathbb
F}(\Omega;C([0,T];\mathbb R^n)) \times L^q_{\mathbb
F}(\Omega;L^2(0,T;\mathbb R^{n\times d}))$ be a solution to
\eqref{Sec3_Haml_Sys} with
\begin{equation}\label{Sec4_tilde_u}
\tilde u(\cdot) = \Gamma(\mathbb K^* \tilde Y(T)) =
\Gamma\Big(B(\cdot)^\top \tilde Y(\cdot) +\sum_{k=1}^d
D_k(\cdot)^\top \tilde Z_k(\cdot)\Big),
\end{equation}
and $\Gamma(\cdot)$ is given by \eqref{G(f)*}. When $p\les
\frac{2\si\rho}{\si\rho-2\rho+2\si}$, Lemma \ref{Lemma 4.6} implies
$\tilde u(\cdot)\in \mathbb U^{p,\rho,\si}[0,T]$. Next we prove that
$\tilde \eta =\tilde Y(T)$ is an optimal solution to Problem (O).
For each $\eta\in L^q_{\mathcal F_T}(\Omega;\mathbb R^n)$, we denote
the solution to BSDE \eqref{Sec3_Dual_Sys} by $(Y(\cdot),Z(\cdot))$,
then
$$
\begin{aligned}
& J(\eta;x,\xi) -J(\tilde\eta;x,\xi) = \frac 1 2 \Big( \| \mathbb
K^*\eta \|^2_{\mathbb U^{p,\rho,\si}[0,T]^*} - \| \mathbb
K^*\tilde\eta \|^2_{\mathbb U^{p,\rho,\si}[0,T]^*} \Big) + \langle
x,\ \mathbb K^*_0\eta-\mathbb K^*_0\tilde \eta \rangle -\mathbb
E\langle \xi,\ \eta-\tilde \eta \rangle.
\end{aligned}
$$
Due to the convexity of $\|\cdot\|^2_{\mathbb
U^{p,\rho,\si}[0,T]^*}$, \eqref{Sec4_tilde_u}, \eqref{K*} and
\eqref{K0*}, we have
$$
\begin{aligned}
& J(\eta;x,\xi) -J(\tilde\eta;x,\xi)\ges  \mathbb E\int_0^T \langle \Gamma(\mathbb
K^*\tilde\eta)(t),\ (\mathbb K^*\eta-\mathbb K^*\tilde \eta)(t)
\rangle dt + \langle x,\ \mathbb K^*_0\eta-\mathbb K^*_0\tilde \eta
\rangle -\mathbb E\langle \xi,\ \eta-\tilde \eta \rangle\\
& =  \mathbb E\int_0^T \Big\langle \tilde u(t),\
B(t)^\top(Y(t)-\tilde Y(t)) +\sum_{k=1}^d D_k(t)^\top (Z_k(t)-\tilde
Z_k(t)) \Big\rangle dt\\
& \qquad +\langle x,\ Y(0)-\tilde Y(0) \rangle -\mathbb E\langle
\xi,\ \eta-\tilde\eta \rangle.
\end{aligned}
$$
We apply It\^{o}'s formula to $\langle \tilde X(\cdot),\
Y(\cdot)-\tilde Y(\cdot) \rangle$, and by \eqref{Sec3_Haml_Sys}, we
obtain the right hand side of the above inequality equals zero.
Hence we have $J(\eta;x,\xi) - J(\tilde \eta;x,\xi) \geq 0$, which
implies that $\tilde \eta = \tilde Y(T)$ is an optimal solution to
Problem (O).

Now, let $(X^i(\cdot), Y^i(\cdot), Z^i(\cdot)) \in L^p_{\mathbb
F}(\Omega;C(0,T;\mathbb R^n)) \times L^q_{\mathbb
F}(\Omega;C(0,T;\mathbb R^n)) \times L^q_{\mathbb
F}(\Omega;L^2(0,T;\mathbb R^{n\times d}))$ ($i=1,2$) be two
solutions to \eqref{Sec3_Haml_Sys}. From the above analysis, both
$Y^1(T)$ and $Y^2(T)$ are optimal controls to Problem (O). By the
uniqueness of optimal control (see Proposition \ref{onto}),
$Y^1(T)=Y^2(T)$. Moreover, by the uniqueness of BSDE, we have
$(Y^1(\cdot), Z^1(\cdot)) = (Y^2(\cdot),Z^2(\cdot))$. Furthermore,
by the uniqueness of FSDE, we have $X^1(\cdot)=X^2(\cdot)$. We
obtain the uniqueness of \eqref{Sec3_Haml_Sys}, and complete the
proof.
\endpf

\ms

\br{Remark 4.8} \rm The notion of adaptability represents a
fundamental difference between deterministic and stochastic systems.
From the derivation of Fr\'echet derivative (see the third line of
\eqref{sec3_Eq_Temp2.5}), we can obtain naturally a process:
$$
\tilde\Gamma(\varphi(\cdot))(t) \equiv \| \varphi(\cdot) \|^{2-\bar
q}_{L^{\bar q}_{\mathbb F}(\Omega;L^\nu(0,T;\mathbb R^m))} \|
\varphi(\cdot,\omega) \|^{\bar q-\nu}_{L^\nu(0,T;\mathbb R^m)}
|\varphi(t)|^{\nu-2}\varphi(t),\quad t\in [0,T]
$$
which is closely linked to our problem. But unfortunately,
$\tilde\Gamma(\varphi(\cdot))(\cdot)$ is not $\mathbb F$-adapted
when $p\neq\frac{2\si\rho}{\si\rho-2\rho+2\si}$ (equivalently $\bar
q\neq \nu$). Hence, in order to meet the requirement of
adaptability, we use
$$
\Gamma(\varphi(\cdot))(t) = \mathbb
E_t[\tilde\Gamma(\varphi(\cdot))(t)],\quad t\in [0,T]
$$
to replace
$\tilde\Gamma(\varphi(\cdot))(\cdot)$. However, this treatment leads
to some difficulty. As a matter of fact, through a direct calculation, we can
obtain that the following equation
$$\left\{ \mathbb E\left( \int_0^T |\tilde\Gamma(\varphi(\cdot))(t)|^{\mu} dt \right)^{\frac {\bar p}
{\mu}} \right\}^{\frac {1} {\bar p}} = \|\varphi(\cdot)\|_{\mathbb
U^{p,\rho,\si}[0,T]^*}$$
holds for any $p\in (1,\infty)$. But due to the introduction of
conditional expectation, we only get an inequality
$$\left\{ \mathbb E\left( \int_0^T |\Gamma(\varphi(\cdot))(t)|^{\mu} dt \right)^{\frac {\bar p}{\mu}}
\right\}^{\frac {1} {\bar p}} \leq \frac{\bar q}{\nu}
\|\varphi(\cdot)\|_{\mathbb U^{p,\rho,\si}[0,T]^*}$$
for $p\les\frac{2\si\rho}{\si\rho-2\rho+2\si}$ (see \eqref{doob}).
The technique involving Doob's martingale inequality used in
Lemma \ref{Lemma 4.6} is invalid for
$p>\frac{2\si\rho}{\si\rho-2\rho+2\si}$. \er

It is naturally to ask what happens if $p >
\frac{2\si\rho}{\si\rho-2\rho+2\si}$? To obtain a similar result as
Theorem \ref{Theorem 4.7}, instead of functional $J(x,\xi;\eta)$
defined by \eqref{Sec3_Dual_Cost}, we introduce the following
functional:
\bel{J'}J'(\eta;x,\xi)={1\over2}\|\dbK^*\eta\|_{\dbU_r^{p,\rho,\si}[0,T]^*}^2+\lan
x,\dbK_0^*\eta\ran-\dbE\lan\xi,\eta \ran,\qq\forall\eta\in
L^q_{\cF_T}(\O;\dbR^n),\ee
where $\dbK^*$ and $\dbK_0^*$ are given by \eqref{K*} and \eqref{K0*}, respectively. Equivalently,
\bel{J'*}\ba{ll}
\ns\ds J'(\eta;x,\xi)={1\over2}\|B(\cd)^\top
Y(\cd)+\sum_{k=1}^dD_k(\cd)^\top
Z_k(\cd)\|_{\dbU_r^{p,\rho,\si}[0,T]^*}^2+\lan x,Y(0)\ran
-\dbE\lan\xi,\eta\ran,\\
\ns\ds\qq\qq\qq\qq\qq\qq\qq\qq\qq\qq\qq\qq\forall\eta\in L^q_{\cF_T}(\O;\dbR^n),\ea\ee
with $(Y(\cd),Z(\cd))$ being the adapted solution to BSDE
\eqref{Sec3_Dual_Sys}. Note that we have changed from
$\dbU^{p,\rho,\si}[0,T]$ to $\dbU_r^{p,\rho,\si}[0,T]$ in the above.
We now pose the following optimization problem.

\ms

\bf Problem (O)$'$. \rm Minimize \eqref{J'} subject to BSDE
\eqref{Sec3_Dual_Sys} over $L^q_{\cF_T}(\O;\dbR^n)$.

\ms

Suppose $\f(\cd), \psi(\cdot)\in L^\n_\dbF(0,T;L^{\bar
q}(\O;\dbR^m))$. A similar procedure as Problem (O) leads to
$$
\begin{aligned}
& \frac 1 2 \frac{d}{d\delta}\Big\{
\|\varphi(\cdot)+\delta\psi(\cdot)\|^2_{L^\nu_{\mathbb
F}(0,T;L^{\bar q}(\Omega;\mathbb R^m))}\Big\}\Big|_{\delta =0} \\
& = \|\varphi(\cdot)\|^{2-\nu}_{L^\nu_{\mathbb F}(0,T;L^{\bar
q}(\Omega;\mathbb R^m))} \int_0^T \Big(\mathbb E |\varphi(t)|^{\bar
q}\Big)^{\frac{\nu-\bar q}{\bar q}} \mathbb E\Big[
|\varphi(t)|^{\bar q-2} \langle \varphi(t),\ \psi(t) \rangle \Big]
dt\\
& \equiv \mathbb E\int_0^T \langle \Gamma'(\varphi(\cdot))(t),\
\psi(t) \rangle dt,
\end{aligned}
$$
where \bel{G'(f)*}\G'(\f(\cd))(t)=
\|\varphi(\cdot)\|^{2-\nu}_{L^\nu_{\mathbb F}(0,T;L^{\bar
q}(\Omega;\mathbb R^m))} \(\dbE|\f(t)|^{\bar q}\)^{{\n-\bar
q\over\bar q}}|\f(t)|^{\bar q-2}\f(t),\qq t\in[0,T]. \ee
Unlike the case of Problem (O), in the above we do not need
conditional expectation. Due to this, the following lemma
holds without the constraint $p\les
\frac{2\si\rho}{\si\rho-2\rho+2\si}$. However, due to the use of
$\dbU^{p,\rho,\si}_r[0,T]^*$, condition $p\ges 2$ is needed.

\bl{Lemma 4.8} \sl Let $p\ges 2$. Then for any $\f(\cd)\in
L^\n_\dbF(0,T;L^{\bar
q}(\O;\dbR^m))\equiv\dbU_r^{p,\rho,\si}[0,T]^*$,
\bel{}\G'(\f(\cd))(\cd)\in L^\m_\dbF(0,T;L^{\bar
p}(\O;\dbR^m))=\dbU_r^{p,\rho,\si}[0,T].\ee

\el

\it Proof. \rm It suffices to calculate the following:
$$\ba{ll}
\ns\ds\int_0^T\Big\{\dbE\[\(\dbE|\f(t)|^{\bar q}\)^{{\n-\bar q\over\bar q}}|\f(t)|^{\bar q-1}\]^{\bar p}\Big\}^{\m\over\bar p}dt=\int_0^T\[\(\dbE|\f(t)|^{\bar q}\)^{\n-\bar q\over\bar q-1}\dbE|\f(t)|^{\bar q}\]^{\m\over\bar p}dt=\int_0^T\(\dbE|\f(t)|^{\bar q}\)^{\n\over\bar q}dt<\infty.\ea$$
Hence, our conclusion follows. \endpf

\ms

Then similar to Theorem \ref{Theorem 4.7}, we have the following result.

\bt{Theorem 4.9} \sl Let {\rm(H1)} and {\rm(H2)$'$} hold and $p\ges 2$. Then
the observability inequality \eqref{weak observability inequality}
holds if and only if for any $(x,\xi)\in\dbR^n\times
L^p_{\cF_T}(\O;\dbR^n)$, Problem {\rm(O)$'$} admits a unique optimal
solution $\bar\eta'\in L^q_{\cF_T}(\O;\dbR^n)$. In this case, the
control $\bar u'(\cd)\in\dbU_r^{p,\rho,\si}[0,T]$ defined by the
following
\bel{baru}\bar u'(\cd)=\G'(\dbK^*\bar\eta'),\ee
with $\G'(\cdot)$ given by \eqref{G'(f)*} steers $x$ to $\xi$.
Moreover, with such defined $\bar u'(\cd)$, the coupled FBSDE
\eqref{Sec3_Haml_Sys} admits a unique adapted solution
$$
(\bar X(\cd),\bar Y(\cd),\bar Z(\cd))\in
L^p_\dbF(\O;C([0,T];\dbR^n))\times
L^q_\dbF(\O;C([0,T];\dbR^n))\times L^q_\dbF(\O;L^2(0,T;\dbR^{n\times
d})).
$$
\et

\section{Norm optimal control problems}\label{S_Norm_Optimal_Control}

When the system \eqref{FSDE1} is $L^p$-exactly controllable on
$[0,T]$ by $\dbU^{p,\rho,\si}[0,T]$ (respectively,
$\dbU_r^{p,\rho,\si}[0,T]$) for any given $p$ satisfying $ p\les
\frac{2\si\rho}{\si\rho-2\rho+2\si}$ (respectively, $p\ges 2$), then
for any $(x, \xi)\in \mathbb R^n\times L^p_{\mathcal
F_T}(\Omega;\mathbb R^n)$, from Theorem \ref{Theorem 4.7}
(respectively, Theorem \ref{Theorem 4.9}), we know that $\bar u$
(respectively, $\bar u'$) defined by \eqref{bar u} (respectively, by
\eqref{baru}) is one of $\mathbb U^{p,\rho,\si}[0,T]$-
(respectively, $\mathbb U_r^{p,\rho,\si}[0,T]$-) admissible controls
which steers the state process from the initial value $x$ to the
terminal value $\xi$. In this section, we shall restrict to the case
$p= \frac{2\si\rho}{\si\rho-2\rho+2\si}$ (respectively, $p\ges 2$)
and further show that $\bar u$ (respectively, $\bar u'$) has a
characteristic of minimum norm. \ms

First, for any given $1<\rho\les\infty$, $2<\si\les\infty$ and $p=
\frac{2\si\rho}{\si\rho-2\rho+2\si}$, we introduce a $\mathbb
U^{p,\rho,\si}[0,T]$-norm optimal control problem: for any
$(x,\xi)\in\mathbb R^n\times L^p_{\mathcal F_T}(\Omega;\mathbb
R^n)$, minimize $\|u\|_{\dbU^{p,\rho,\si}[0,T]}$ over the $\mathbb
U^{p,\rho,\si}[0,T]$-admissible control set:
$$
\mathcal U(x,\xi) := \Big\{ u(\cdot)\in \dbU^{p,\rho,\si}[0,T]\
\Big|\ X(T;x,u(\cdot)) = \xi \Big\}.
$$
For simplicity, we denote the above $\mathbb
U^{p,\rho,\si}[0,T]$-norm optimal control problem by {\bf Problem
(N)}. Note that the system \eqref{FSDE1} is $L^p$-exactly
controllable  on $[0,T]$ by $\dbU^{p,\rho,\si}[0,T]$ if and only if,
for any $(x,\xi)\in\mathbb R^n\times L^p_{\mathcal
F_T}(\Omega;\mathbb R^n)$, the $\mathbb
U^{p,\rho,\si}[0,T]$-admissible control set $\mathcal U(x,\xi)$ is
not empty. We call $\bar u\in \mathcal U(x,\xi)$ a $\mathbb
U^{p,\rho,\si}[0,T]$-norm optimal control to  Problem (N) if
$$\|\bar u(\cdot)\|_{\dbU^{p,\rho,\si}[0,T]} = \inf_{u(\cdot)\in \mathcal U(x,\xi)}
\|u(\cdot)\|_{\dbU^{p,\rho,\si}[0,T]}.$$

For any given $2\les p <\rho\les \infty$ and $2< \si\les \infty$, we
similarly introduce the $\mathbb U_r^{p,\rho,\si}[0,T]$-norm optimal
control problem reading:

\ms

{\bf Problem (N)$'$.} For any $(x,\xi)\in\mathbb R^n\times
L^p_{\mathcal F_T}(\Omega;\mathbb R^n)$, minimize
$\|\cdot\|_{\dbU_r^{p,\rho,\si}[0,T]}$ over the $\mathbb
U_r^{p,\rho,\si}[0,T]$-admissible control set
$$
\mathcal U'(x,\xi) := \Big\{ u\in \dbU_r^{p,\rho,\si}[0,T]\ \Big|\
X(T;x,u(\cdot)) = \xi \Big\}.
$$

In the previous section, we have given some equivalent conditions
for the $L^p$-exact controllability of system \eqref{FSDE1} on
$[0,T]$ by $\mathbb U^{p,\rho,\si}[0,T]$ (respectively, $\mathbb
U^{p,\rho,\si}_r[0,T]$) (see Theorem \ref{equivalence}, Theorem
\ref{observe}, Theorem \ref{Theorem 4.7} and Theorem \ref{Theorem
4.9}). Now, by virtue of the related $\mathbb U^{p,\rho,\si}[0,T]$-
(respectively, $\mathbb U_r^{p,\rho,\si}[0,T]$-) norm optimal
control problem, we present another one.

\bt{Sec4_Theorem_Equiv} \sl Let the assumptions {\rm(H1)}, {\rm(H2)}
and $p= \frac{2\si\rho}{\si\rho-2\rho+2\si}$ hold. Then the
following two statements are equivalent:

\ms

$\bullet$ For any $(x,\xi)\in \mathbb R^n\times L^p_{\mathcal
F_T}(\Omega;\mathbb R^n)$, Problem {\rm(O)}  admits a unique optimal
solution $\bar\eta\in L^q_{\mathcal F_T}(\O;\mathbb R^n)$;

\ms

$\bullet$ For any $(x,\xi)\in \mathbb R^n\times L^p_{\mathcal
F_T}(\Omega;\mathbb R^n)$, Problem {\rm(N)} admits a unique optimal
control $\bar u(\cdot)\in \dbU^{p,\rho,\si}[0,T]$.

\noindent Moreover, the unique norm optimal control $\bar u$ to
Problem {\rm(N)} is given by \eqref{bar u},
and the minimal norm  is given by
\begin{equation}\label{Sec4_Min_Norm}
\|\bar u(\cdot)\|_{\dbU^{p,\rho,\si}[0,T]} = \sqrt{\mathbb E\langle
\xi,\ \bar\eta\rangle -\langle x,\ \mathbb K_0^*\bar\eta \rangle}.
\end{equation}
The minimal value of functional $J(\cdot;x,\xi)$ is given by
\begin{equation}\label{Sec4_Min_J}
J(\bar\eta;x,\xi) = -\frac 1 2 \|\bar
u(\cdot)\|^2_{\dbU^{p,\rho,\si}[0,T]} = -\frac 1 2 \Big(\mathbb
E\langle \xi,\ \bar\eta\rangle -\langle x,\ \mathbb K_0^*\bar\eta
\rangle\Big).
\end{equation}
\et

\it Proof. \rm (Sufficiency). Since for any $(x,\xi)\in\mathbb
R^n\times L^p_{\mathcal F_T}(\Omega;\mathbb R^n)$, Problem (N)
admits an optimal control, then the corresponding $\mathbb
U^{p,\rho,\si}[0,T]$-admissible control set $\mathcal U(x,\xi)$ is
not empty. Therefore, the system \eqref{FSDE1} is $L^p$-exactly
controllable. By Theorem \ref{equivalence} and Theorem \ref{Theorem
4.7}, the first statement holds true.

(Necessity). First of all, when $p=
\frac{2\si\rho}{\si\rho-2\rho+2\si}$, by Lemma \ref{Lemma 4.6},
$\bar u$ defined by \eqref{bar u} is a $\mathbb
U^{p,\rho,\si}[0,T]$-admissible control to Problem (N). Let
$\bar\eta\in L^q_{\mathcal F_T}(\Omega;\mathbb R^n)$ be the optimal
solution to Problem (O), and $(\bar Y(\cdot),\bar Z(\cdot))$ be the
solution to BSDE \eqref{Sec3_Dual_Sys} with terminal condition $\bar
Y(T)=\bar\eta$. For any $\eta\in L^q_{\mathcal F_T}(\Omega;\mathbb
R^n)$, similarly, we denote by $(Y(\cdot),Z(\cdot))$ the solution to
BSDE \eqref{Sec3_Dual_Sys} with terminal condition $Y(T)=\eta$.
Since $\bar\eta$ is optimal, we obtain the following relation named
the Euler-Lagrange equation:
\begin{equation}\label{Sec3_EulerLagrange}
\begin{aligned}
0=\ & \frac{d}{d\delta}
J(\bar\eta+\delta\eta;x,\xi)\Big|_{\delta =0}=\dbE\int_0^T \langle \Gamma(\mathbb K^*\bar\eta)(t),\
(\mathbb K^*\eta)(t) \rangle dt +\langle x,\ \mathbb K^*_0 \eta
\rangle
-\mathbb E\langle \xi,\ \eta \rangle\\
=\ & \mathbb E\int_0^T \langle \bar u(t),\ (\mathbb K^*\eta)(t)
\rangle dt +\langle x, Y(0) \rangle -\mathbb E\langle \xi,\ \eta
\rangle.
\end{aligned}
\end{equation}
Meanwhile, for any control $u(\cdot)\in \mathcal U(x,\xi)$, by applying
It\^o's formula to $\lan X(\cd\,;x,u(\cd)),Y(\cd)\ran$ on the interval $[0,T]$, we have
\begin{equation}\label{ito2}
\mathbb E\int_0^T \langle  u(t),\ (\mathbb K^*\eta)(t) \rangle dt
+\langle x,\ Y(0) \rangle -\mathbb E\langle \xi,\ \eta \rangle =0.
\end{equation}
By
letting $\eta = \bar\eta$ in both \eqref{Sec3_EulerLagrange} and
\eqref{ito2}, for any $u\in\mathcal U(x,\xi)$, we obtain
\begin{equation}\label{Sec4_Temp1}
\mathbb E \int_0^T \lan \bar u(t),\ (\mathbb K^*\bar\eta)(t)\ran dt
= \mathbb E\left\langle \xi,\ \bar\eta \right\rangle -\langle x,\
\bar Y(0) \rangle = \mathbb E\int_0^T \langle u(t),\ (\mathbb
K^*\bar\eta)(t) \rangle dt.
\end{equation}
It is easy to calculate that
$$
\mathbb E \int_0^T \lan \bar u(t), (\mathbb K^*\bar\eta)(t)\ran
dt=\|\mathbb K^*\bar\eta\|_{\mathbb U^{p,\rho,\si}[0,T]^*}^2,
$$
and
$$
\|\bar u(\cdot)\|_{\dbU^{p,\rho,\si}[0,T]}=\|\mathbb
K^*\bar\eta\|_{\mathbb U^{p,\rho,\si}[0,T]^*}.
$$
By above three equations and
H\"{o}lder's inequality,
we get the optimality of $\bar u$. The uniqueness of the optimal
control to Problem (N) comes from the strictly convex property  of
$\dbU^{p,\rho,\si}[0,T]$.


Let us come back to the Euler-Lagrange equation. Letting
$\eta=\bar\eta$ and from the definition
 of $\bar u$,
\eqref{Sec3_EulerLagrange} is reduced to
$$\|\bar u(\cdot)\|_{\dbU^{p,\rho,\si}[0,T]}^2 = \mathbb E\langle \xi,\
\bar\eta\rangle -\langle x,\ \bar Y(0) \rangle,$$
which, together with \eqref{K0*},  implies \eqref{Sec4_Min_Norm}. Then, we calculate
$$J(\bar\eta;x,\xi) = \frac 1 2 \|\bar u(\cdot)\|_{\dbU^{p,\rho,\si}[0,T]}^2
-\Big( \mathbb E\langle \xi,\ \bar\eta\rangle -\langle x,\ \bar Y(0)
\rangle \Big) = - \frac 1 2 \|\bar
u(\cdot)\|_{\dbU^{p,\rho,\si}[0,T]}^2,$$
which is \eqref{Sec4_Min_J}, and the proof is completed.
\endpf

\br{Remark 5.2} \rm When we consider the corresponding $\mathbb
U^{p,\rho,\si}[0,T]$-norm optimal control problems with $p<
\frac{2\si\rho}{\si\rho-2\rho+2\si}$ in the same way, we can obtain
$$\|\mathbb K^*\bar\eta\|_{\mathbb U^{p,\rho,\si}[0,T]^*} \les
\|u(\cdot)\|_{\dbU^{p,\rho,\si}[0,T]},$$
where $u(\cdot)$ is any given $\mathbb
U^{p,\rho,\si}[0,T]$-admissible control. If we can also obtain
$$\|\bar u(\cdot)\|_{\dbU^{p,\rho,\si}[0,T]} \les
\|\mathbb K^*\bar\eta\|_{\mathbb U^{p,\rho,\si}[0,T]^*},$$
then we will solve $\mathbb U^{p,\rho,\si}[0,T]$-norm optimal
control problems. However, due to the additional conditional
expectation in the Fr\'echet derivative $\Gamma(\cdot)$, we cannot
prove the above inequality. In fact, we only have (comparing \eqref{doob})
$$\|\bar u(\cdot)\|_{\mathbb U^{p,\rho,\si}[0,T]} \les \frac{\bar q}{\nu}
\|\mathbb K^*\bar\eta\|_{\mathbb U^{p,\rho,\si}[0,T]^*}.$$
\er

Since conditional expectation is not introduced in $\Gamma'(\cdot)$,
we can solve the $\mathbb U_r^{p,\rho,\si}[0,T]$-norm optimal
control problems for any $p\ges 2$, whose proof is similar to
that of Theorem \ref{Sec4_Theorem_Equiv}.

\bt{Sec4_weak_Theorem_Equiv} \sl Let the assumptions {\rm(H1)},
{\rm(H2)$'$} and $p\ges 2$ hold. Then the following two statements
are equivalent:


\ms

$\bullet$ For any $(x,\xi)\in \mathbb R^n\times L^p_{\mathcal
F_T}(\Omega;\mathbb R^n)$, Problem {\rm(O)$'$}  admits a unique optimal solution
$\bar\eta'\in L^q_{\mathcal F_T}(\O;\mathbb R^n)$;

\ms

$\bullet$ For any $(x,\xi)\in \mathbb R^n\times L^p_{\mathcal
F_T}(\Omega;\mathbb R^n)$, Problem {\rm(N)$'$} admits a unique
optimal control $\bar u'(\cdot)\in \dbU_r^{p,\rho,\si}[0,T]$.

\noindent Moreover, the unique norm optimal control $\bar u'(\cdot)$
to Problem {\rm(N)$'$} is given by \eqref{baru},
and the minimal norm  is given by
\begin{equation}
\|\bar u'(\cdot)\|_{\dbU_r^{p,\rho,\si}[0,T]} = \sqrt{\mathbb
E\langle \xi,\ \bar\eta'\rangle -\langle x,\ \mathbb K_0^*\bar\eta'
\rangle}.
\end{equation}
The minimal value of functional $J'(\cdot;x,\xi)$ is given by
\begin{equation}
J'(\bar\eta;x,\xi) = -\frac 1 2 \|\bar
u'(\cdot)\|^2_{\dbU_r^{p,\rho,\si}[0,T]} = -\frac 1 2 \Big(\mathbb
E\langle \xi,\ \bar\eta'\rangle -\langle x,\ \mathbb K_0^*\bar\eta'
\rangle\Big).
\end{equation}
\et

\br{Remark 5.4} \rm When $p=\frac{2\si\rho}{\si\rho-2\rho+2\si}$
(equivalently, $\bar p =\mu$), the spaces $\mathbb
U^{p,\rho,\si}[0,T]$ and $\mathbb U^{p,\rho,\si}_r[0,T]$ coincide with $L^\mu_{\mathbb F}(0,T;\mathbb R^m)$, and then both
Problem (N) and Problem (N)' imply the $L^\mu_{\mathbb
F}(0,T;\mathbb R^m)$-norm optimal control problem. Precisely,
Theorems \ref{Sec4_Theorem_Equiv} and \ref{Sec4_weak_Theorem_Equiv} provide the same result for the
$L^\mu_{\mathbb F}(0,T;\mathbb R^m)$-norm optimal control problem
when the system is $L^p$-exactly controllable with $p\ges 2$.
However, when $1<p<2$, only Theorem \ref{Sec4_Theorem_Equiv} gives a
result, and Theorem \ref{Sec4_weak_Theorem_Equiv} is invalid.
\er



More specifically, when $\si=\infty$ and $p=\frac{2\rho}{\rho+2}$
(equivalently, $\bar p =\mu =2$), $\mathbb U^{p,\rho,\si}[0,T]
=\mathbb U^{p,\rho,\si}_r[0,T] = L^2_{\mathbb F}(0,T;\mathbb R^m)$.
Problem (N) becomes the classical norm optimal control problem (see
\cite{Wang-Zhang 2015}), and Theorem \ref{Sec4_Theorem_Equiv}
provides a result for the classical $L^2_{\mathbb F}(0,T;\mathbb
R^m)$-norm optimal control problem. We notice that in the present
paper the matrices $B(\cdot)$ and $D_k(\cdot)$ ($k=1,2\dots,d$) are
not necessary to be bounded (see Assumption (H2)), while in the
literature, only bounded matrix cases were studied. Furthermore, instead
of the standard norm $\|\cd\|_{L_{\mathbb F}^{2}(0,T;\mathbb R^m)}$,
we can extend our method to minimize the following generalized
weighted norm
\begin{equation}\label{Sec4_Weighted_Norm}
\left(\mathbb E\int_0^T \langle R(t)u(t),u(t) \rangle
dt\right)^{\frac 1 2}
\end{equation}
subject to \eqref{FSDE1} over the $L^2_{\mathbb F}(0,T;\mathbb
R^m)$-admissible control set $\mathcal U(x,\xi) $, where $R(\cd)\in
L^\infty_{\mathbb F}(0,T;\mathbb R^{m\times m})$ is symmetric, and
there exists a constant $\delta>0$ such that $R(t)-\delta I$ is
positive semi-definite for a.e. $t\in [0,T]$.

In fact, due to the definition of $R(\cd)$ and Denman-Beavers iteration
\cite{Lee-Markus 1967} for the square roots of matrices, there
exists a matrix-valued process $N(\cd)\in L^\infty_{\mathbb
F}(0,T;\mathbb R^{m\times m})$ which is invertible and $N^{-1}$ is
bounded also, such that $R(\cdot) = N^\top(\cdot)N(\cdot)$. Then,
letting
$$\begin{aligned}
& \hat u(\cdot) = N(\cdot)u(\cdot),\q \hat B(\cdot) = B(\cdot)N^{-1}(\cdot),\quad \hat
D_j(\cdot)=D_j(\cdot)N^{-1}(\cdot),\quad j=1,2,\dots,d,
\end{aligned}$$
we transform the original weighted norm optimal control problem
\eqref{FSDE1} and \eqref{Sec4_Weighted_Norm} to the following
equivalent standard one: to minimize $\|\hat
u(\cdot)\|_{L^{2}_{\mathbb F}(0,T;\mathbb R^m)}$ subject to
$$\left\{
\begin{aligned}
d X(t) =\ & \Big[A(t) X(t) +\hat B(t)\hat u(t)\Big] dt +\sum_{k=1}^d
\Big[ C_k(t) X(t) +\hat D_k(t)\hat u(t)
\Big]dW_k(t),\quad t\in [0,T],\\
 X(0) =\ & x,
\end{aligned}
\right.$$
over the corresponding $L^2_{\mathbb F}(0,T;\mathbb R^m)$-admissible
control set
$$
\begin{aligned}
\hat{\mathcal U}(x,\xi) \equiv\ & \Big\{ \hat u(\cdot)\in
L^2_{\mathbb F}(0,T;\mathbb R^m)\ \Big|\  X(T;x,u(\cdot)) = \xi
\Big\}=N^{-1}\mathcal U(x,\xi)\\
\equiv\ & \Big\{ \hat u(\cdot) = N^{-1} u(\cdot)\ \Big|\ u(\cdot)\in
\mathcal U(x,\xi) \Big\}.
\end{aligned}$$
Applying the results of Theorems \ref{equivalence},
\ref{Theorem 4.7} and \ref{Sec4_Theorem_Equiv}, we
can solve the generalized weighted norm optimal control problem.

\bc{Corollary 5.5} \sl Let $\si=\infty$, $p=\frac{2\rho}{\rho+2}$
and {\rm (H1)--(H2)} hold. Then the system
\eqref{FSDE1} is $L^p$-exactly controllable on $[0,T]$ by
$L^2_{\mathbb F}(0,T;\mathbb R^m)$, if and only if, for any $(x,
\xi)\in\mathbb R^n\times L^p_{\mathcal F_T}(\Omega;\mathbb R^n)$,
the weighted norm optimal control problem \eqref{FSDE1} and
\eqref{Sec4_Weighted_Norm} admits a unique optimal control. In this
case, with
\begin{equation}\label{Sec4_Norm_Optimal_Control_Weighted}
\bar u(t) = R^{-1}(t) \Big( B^\top(t)\bar Y(t) +\sum_{k=1}^d
D_k(t)^\top \bar Z_k(t) \Big),\quad t\in [0,T],
\end{equation}
and the following coupled FBSDE
\begin{equation}\label{Sec4_Haml_Sys_Weighted}
\left\{
\begin{aligned}
d\bar X(t) =\ & \Big[ A(t)\bar X(t) +B(t)\bar u(t) \Big] dt
+\sum_{k=1}^d
\Big[ C_k(t)\bar X(t) +D_k(t)\bar u(t) \Big] dW_k(t),\\
-d\bar Y(t) =\ & \Big[ A(t)^\top \bar Y(t) +\sum_{k=1}^d
C_k(t)^\top \bar Z_k(t) \Big] dt -\sum_{k=1}^d \bar Z_k(t) dW_k(t),\\
\bar X(0) =\ & x,\qquad \bar X(T) = \xi
\end{aligned}
\right.
\end{equation}
admits a unique adapted solution
$$
(\bar X(\cdot),\bar Y(\cdot),\bar Z(\cdot)) \in L^p_{\mathbb
F}(\Omega;C(0,T;\mathbb R^n)) \times L^p_{\mathbb
F}(\Omega;C(0,T;\mathbb R^n)) \times L^p_{\mathbb
F}(\Omega;L^2(0,T;\mathbb R^{n\times d})).
$$
Moreover, $\bar u(\cdot)$ defined by
\eqref{Sec4_Norm_Optimal_Control_Weighted} is the unique weighted
norm optimal control, and the minimal weighted norm is given by
$$
\left(\mathbb E\int_0^T \langle R(t)\bar u(t),\bar u(t) \rangle
dt\right)^{\frac 1 2} = \sqrt{\mathbb E \langle \xi,\bar Y(T)
\rangle -\langle x,\bar Y(0) \rangle}.
$$
\ec

\section*{Acknowledgement}

This work was carried out during the stay of Yanqing Wang, Donghui
Yang and Zhiyong Yu at University of Central Florida, USA. They
would like to thank the Department of Mathematics for its
hospitality, and the financial support from China Scholarship Council. Yanqing Wang also gratefully acknowledges Professor Qi
L\"{u}, Sichuan University, for stimulating discussion during this
work.

\end{document}